\documentclass[a4paper, twoside, 11pt]{scrartcl}
\pdfoutput=1

\usepackage{etex}
\usepackage[a4paper, top=88pt, bottom=88pt, left=72pt, right=72pt, headsep=16pt, footskip=28pt]{geometry}
\usepackage{amsfonts,amssymb,amsmath,amsthm,amstext,amssymb,amsopn,mathtools,nicefrac,xfrac}
\usepackage[authoryear,sort,round]{natbib} 
\usepackage[utf8]{inputenc}                 
\usepackage{booktabs}                       
\usepackage{nicefrac}                       
\usepackage{microtype}                      
\usepackage{graphics,graphicx}              
\usepackage{subcaption}                     
\usepackage[ruled,linesnumbered]{algorithm2e} 
\usepackage{verbatim}                       
\usepackage{appendix}                       
\usepackage{bm}                             
\usepackage{array}                          
\usepackage[usenames,dvipsnames]{xcolor}    
\usepackage{hyperref}                       
\hypersetup{colorlinks=true,linkcolor=blue,filecolor=magenta,urlcolor=blue,citecolor=OliveGreen}
\usepackage[noabbrev,capitalise,nosort,nameinlink]{cleveref}  

\newcommand*{\arXiv}[1]{\bgroup\color{blue}\href{https://arxiv.org/abs/#1}{arXiv:#1}\egroup}
\newcommand*{\doi}[1]{\bgroup\color{blue}\href{https://doi.org/#1}{doi:#1}\egroup}
\newcommand*{\email}[1]{\bgroup\color{blue}\href{mailto:#1}{#1}\egroup}
\renewcommand*{\url}[1]{\bgroup\color{blue}\href{#1}{#1}\egroup}
\usepackage{enumitem, moreenum}
\setlist[enumerate]{nosep}
\setlist[itemize]{nosep}
\usepackage{mleftright} \mleftright
\renewcommand{\qedsymbol}{$\blacksquare$}
\renewenvironment{proof}[1][\proofname]{\noindent{\bfseries\sffamily #1.} }{\hfill\qedsymbol\medskip}
\usepackage[labelfont={sf,bf}]{caption}
\DeclareCaptionLabelSeparator{figlabelsep}{\,\,\,}
\usepackage{scrlayer-scrpage, xhfill}
\automark[section]{section}
\setkomafont{pageheadfoot}{\normalcolor\sffamily}
\setkomafont{pagenumber}{\normalfont\normalsize\sffamily}
\clearpairofpagestyles
\let\oldtitle\title
\renewcommand{\title}[1]{\oldtitle{#1}\newcommand{\theshorttitle}{#1}}
\newcommand{\shorttitle}[1]{\renewcommand{\theshorttitle}{#1}}
\let\oldauthor\author
\renewcommand{\author}[1]{\oldauthor{#1}\newcommand{\theshortauthor}{#1}}
\newcommand{\shortauthor}[1]{\renewcommand{\theshortauthor}{#1}}
\cohead{\xrfill[0.525ex]{0.6pt}~\theshorttitle~\xrfill[0.525ex]{0.6pt}}
\cehead{\xrfill[0.525ex]{0.6pt}~\theshortauthor~\xrfill[0.525ex]{0.6pt}}
\newcommand{\theabstract}[1]{\par\bgroup\noindent\textbf{\textsf{Abstract.}} #1\egroup}
\newcommand{\thekeywords}[1]{\par\smallskip\bgroup\noindent\textbf{\textsf{Keywords.}}\newcommand{\and}{ $\bullet$ } #1\egroup}
\newcommand{\themsc}[1]{\par\smallskip\bgroup\noindent\textbf{\textsf{2020 Mathematics Subject Classification.}}\newcommand{\and}{ $\bullet$ } #1\egroup}
\newcommand*{\affilref}[1]{\ref{affiliation#1}}
\newcommand*{\affiliation}[3]{
	\footnotetext[#1]{\label{affiliation#2} #3}
}

\numberwithin{equation}{section}
\numberwithin{figure}{section}
\numberwithin{table}{section}

\newcommand*{\defeq}{\coloneqq}

\renewcommand*{\geq}{\geqslant}

\renewcommand*{\leq}{\leqslant}
\newcommand*{\rd}{\mathrm{d}}
\newcommand*{\Reals}{\mathbb{R}}
\newcommand*{\posE}{\text{E}}
\newcommand*{\negE}{\text{E}\mathord{-}}
\newcommand*{\F}{\mathcal{F}_{\Delta T}}
\newcommand*{\G}{\mathcal{G}_{\Delta T}}

\newtheorem{theorem}{\sffamily Theorem}[section]

\theoremstyle{definition}

\newtheorem{assumption}[theorem]{\sffamily Assumption}

\renewenvironment{proof}[1][\proofname]{\noindent{\bfseries\sffamily #1.} }{\hfill\qedsymbol\medskip}

\crefname{assumption}{Assumption}{Assumptions}
\Crefname{assumption}{Assumption}{Assumptions}

\title{GParareal: A time-parallel ODE solver using Gaussian process emulation}
\shorttitle{GParareal}
\author{
    Kamran Pentland\textsuperscript{\affilref{Warwick}} 
    \and
    Massimiliano Tamborrino\textsuperscript{\affilref{Warwick}} 
    \and
    T.~J.~Sullivan\textsuperscript{\affilref{Warwick},\affilref{Turing}} \\
    \and
    James Buchanan\textsuperscript{\affilref{Culham}} 
    \and
    L.~C.~Appel\textsuperscript{\affilref{Culham}}
}
\shortauthor{K.~Pentland, M.~Tamborrino, T.~J.~Sullivan, J.~Buchanan, and L.~C.~Appel}
\date{\today}

\begin{document}
\maketitle

\affiliation{1}{Warwick}{University of Warwick, Coventry, CV4 7AL, United Kingdom\newline (\email{kamran.pentland@warwick.ac.uk}, \email{massimiliano.tamborrino@warwick.ac.uk}, \email{t.j.sullivan@warwick.ac.uk})}
\affiliation{2}{Turing}{Alan Turing Institute, British Library, 96 Euston Road, London NW1 2DB, United Kingdom}
\affiliation{3}{Culham}{Culham Centre for Fusion Energy, Culham Science Centre, Abingdon, Oxfordshire, OX14 3DB, United Kingdom\newline (\email{james.buchanan@ukaea.uk}, \email{lynton.appel@ukaea.uk})}


\begin{abstract}\small
    \theabstract{%
        Sequential numerical methods for integrating initial value problems (IVPs) can be prohibitively expensive when high numerical accuracy is required over the entire interval of integration. 
        One remedy is to integrate in a parallel fashion, ``predicting'' the solution serially using a cheap (coarse) solver and ``correcting'' these values using an expensive (fine) solver that runs in parallel on a number of temporal subintervals.
        In this work, we propose a time-parallel algorithm (\emph{GParareal}) that solves IVPs by modelling the correction term, i.e.\ the difference between fine and coarse solutions, using a Gaussian process emulator. 
        This approach compares favourably with the classic \emph{parareal} algorithm and we demonstrate, on a number of IVPs, that GParareal can converge in fewer iterations than parareal, leading to an increase in parallel speed-up. 
        GParareal also manages to locate solutions to certain IVPs where parareal fails and has the additional advantage of being able to use archives of legacy solutions, e.g. solutions from prior runs of the IVP for different initial conditions, to further accelerate convergence of the method --- something that existing time-parallel methods do not do.
    }
    \thekeywords{%
        {Gaussian process emulation}%
        \and%
        {parallel-in-time}%
        \and%
        {parareal}%
        \and%
        {initial value problems}%
    }
    \themsc{
        {65L05}
        \and%
        {65Y05}
        \and%
        {60G15}
    }
\end{abstract}



\section{Introduction} \label{sec:intro}

\subsection{Motivation and background}
This paper is concerned with the numerical solution of a system of $d \in \mathbb{N}$ ordinary differential equations (ODEs) of the form
\begin{equation} \label{eq:ODE}
    \frac{\rd \bm{u}}{\rd t} = \bm{f} \bigl( t,\bm{u}(t) \bigr) \quad \text{over} \quad t \in [t_0, T], \quad \text{with} \quad \bm{u}(t_0) = \bm{u}_0 \in \mathcal{U} \subset \Reals^{d},
\end{equation}
where $\bm{f}\colon [t_0,T] \times \mathcal{U} \to \Reals^d$ is a nonlinear function with sufficiently many continuous partial derivatives, $\bm{u}\colon [t_0, T] \to \mathcal{U}$ is the time-dependent solution, and $\bm{u}_0$ is the initial value at time $t_0$.
We seek numerical solutions $\bm{U}_j \approx \bm{u}(t_j)$ to the initial value problem (IVP) in \eqref{eq:ODE} on a pre-defined mesh $\bm{t} = (t_0,\dots,t_J)$, where $t_{j+1} = t_j + \Delta T$ for fixed $\Delta T = (T-t_0)/J$.

More specifically, we are concerned with IVPs where: (i) the interval of integration, $[t_0,T]$; (ii) the number of mesh points, $J+1$; or (iii) the wallclock time to evaluate the vector field, $\bm{f}$, is so large that such numerical solutions take hours, days, or even weeks to obtain using classical sequential integration methods, e.g.\ implicit/explicit Runge--Kutta methods \citep{hairer1993}.
Expensive vector fields $\bm{f}$ can, for example, arise when (spatially) discretising partial differential equations (PDEs) into a system of ODEs.
Runtime issues also arise when solving IVPs with spatial or other non-temporal dependencies in that, even though highly efficient domain decomposition methods exist \citep{dolean2015}, the parallel speed-up of such methods on high performance computers (HPCs) is still constrained by the serial nature of the time-stepping scheme.
Therefore, with the advent of exascale HPCs on the horizon \citep{mann2020}, there has been renewed interest in developing more efficient and robust \emph{time-parallel} algorithms to reduce wallclock runtimes for IVP simulations in applications spanning numerical weather prediction \citep{hamon2020}, kinematic dynamo modelling \citep{clarke2020}, and plasma physics \citep{samaddar2010,samaddar2019} to name but a few.
In this work, we focus on the development of such a time-parallel method. 

To solve \eqref{eq:ODE} in parallel, one must overcome the causality principle of time: solutions at later times depend on solutions at earlier times.
In recent years, a growing number of time-parallel algorithms, whereby one partitions $[t_0,T]$ into $J$ `slices' and attempts to solve $J$ smaller IVPs using $J$ processors, have been developed to speed-up IVP simulations;
see \citet{gander2015} and \citet{ong2020} for comprehensive reviews.
We take inspiration from the \emph{parareal} algorithm \citep{lions2001}, a multiple shooting-type (or multigrid \citep{gander2007}) method that uses a predictor-corrector update rule based on two numerical integrators, one coarse- and one fine-grained in time, to iteratively locate solutions $\bm{U}^k_j$ to \eqref{eq:ODE} in parallel.
At any iteration $k \in \{1,\dots,J\}$ of parareal, the `correction' is given by the residual between fine and coarse solutions obtained during iteration $k-1$ (further details are provided in \Cref{sec:parareal}).
In a Markovian-like manner, all fine/coarse information about the solution obtained prior to iteration $k-1$ is ignored by the predictor-corrector rule, a feature present in most parareal-type algorithms and variants \citep{elwasif2011,ait-ameur2020,maday2020,dai2013,pentland2022}.
Our goal is to demonstrate that such ``acquisition'' data, i.e.\ fine and coarse solution information accumulated up to iteration $k$, can be exploited using a statistical \emph{emulator} in order to determine a solution in faster wallclock time than parareal.
Making use of acquisition data in parareal is mentioned briefly in the appendix of \cite{maday2020}, in the context of spatial domain decomposition and high-order time-stepping, but has yet to be investigated further.


In particular, we use a Gaussian process (GP) emulator \citep{ohagan1978,rasmussen2004} to rapidly infer the (expensive-to-simulate) multi-fidelity correction term in parareal.
The emulator is trained using acquisition data from \emph{all} prior iterations, with data from the fine integrator having been obtained in parallel.
Similar to parareal, we derive a predictor-corrector-type scheme where the coarse integrator makes rapid low-accuracy predictions about the solutions which are subsequently refined using a correction, now inferred from the GP emulator.
In addition to using an emulator, the difference between this approach and parareal is that the new correction term is formed of integrated solutions values at the current iteration $k$, rather than $k-1$.
Supposing that the fine solver is of sufficient accuracy to exactly solve the IVP, the algorithm presented in this paper determines numerical solutions $\bm{U}^k_j$ that converge (assuming the emulator is sufficiently well trained) toward the exact solutions $\bm{U}_j$ over a number of iterations.
This new approach is particularly beneficial if one wishes to fully understand and evaluate the dynamics of \eqref{eq:ODE} by simulating solutions for a range of initial values $\bm{u}_0$ or over different time intervals.
Firstly, if one can obtain additional parallel speedup, generating such a sequence of independent simulations becomes more computationally tractable in feasible time. 
Secondly, the ``legacy'' data, i.e.\ solution information accumulated between independent simulations, can be used to inform future simulations by increasing the size of the dataset available to the emulator.
Being able to re-use (expensive) acquisition or legacy data to integrate IVPs such as \eqref{eq:ODE} in parallel is not something, to the best of our knowledge, that existing time-parallel algorithms currently do.


In recent years, there has been a surge in interest in the field of \emph{probabilistic numerics} \citep{hennig2022,oates2019}, where ``ODE filters'' have been developed to solve ODEs using GP regression techniques.
Instead of calculating a numerical solution on the mesh $\bm{t}$, as classical integration methods do, ODE filters return a probability measure over the solution at any $t \in [t_0,T]$ \citep{schober2019,tronarp2019,bosch2021,wenger2021}.
Such methods solve sequentially in time, conditioning the GP on acquisition data, i.e.\ solution and derivative evaluations, at competitive computational cost (compared to classical methods) \citep{kersting2020,kramer2021}.
However, integrating IVPs with large time intervals or expensive vector fields using such filters is still a computationally intractable process. 
As such, our aim is to fuse aspects of time-parallelism with the Bayesian methods showcased in ODE filters---something briefly mentioned in \cite{kersting2018} and \cite{pentland2022}, but not yet explored. 
Whereas ODE filters use GPs to explicitly model the \emph{solution} to an IVP, we instead use them to model the \emph{residual} between approximate solutions provided by the deterministic fine and coarse solvers, i.e.\ the parareal correction.
While the method proposed in this paper \emph{does not} return a probabilistic solution to \eqref{eq:ODE}, we believe that it constitutes a positive step in this direction.


\subsection{Contributions and outline}
The rest of this paper is structured as follows.
In \Cref{sec:parareal}, we introduce parareal, providing an overview of the algorithm and its computational complexity for a scalar ODE.
In \Cref{sec:GParareal}, we present our algorithm, henceforth referred to as GParareal, in which we describe how a GP emulator, conditioned on acquisition data obtained in parallel throughout the simulation, is used to refine coarse numerical solutions to a scalar ODE.
In addition, we detail the computational complexity of GParareal, provide a bound for its numerical error at a given iteration, and describe the extension for solving systems of ODEs.
Numerical experiments are performed using HPC facilities in \Cref{sec:numerics}. 
We demonstrate good performance of GParareal against parareal in terms of convergence, wallclock time, and solution accuracy on a number of low-dimensional ODE problems using just acquisition data.
Furthermore, we demonstrate how the GP emulator enables convergence in cases where the coarse solver is too inaccurate for parareal to converge and show that legacy simulation data can be used to obtain solutions even faster, retaining comparable numerical accuracy.
We discuss the benefits, drawbacks, and open questions surrounding GParareal in \Cref{sec:discussion}.


\section{Parareal} \label{sec:parareal}

Here we briefly recall the parareal algorithm \citep{lions2001}, first describing the fine- and coarse-grained numerical solvers it uses, then the algorithm itself, and finally some remarks on complexity, numerical speed-up, and choice of solvers. 
For a full mathematical derivation and exposition of parareal, refer to \cite{gander2007}. To simplify notation, we describe parareal for solving a scalar-valued autonomous ODE, i.e.\ $\bm{f}(t,\bm{u}(t)) \defeq f(u(t))$ in \eqref{eq:ODE}, without loss of generality.

\subsection{The solvers}
To calculate a solution to \eqref{eq:ODE}, parareal uses two one-step\footnote{Multi-step numerical integrators have been tested with parareal \citep{ait-ameur2020,ait-ameur2021}.
However, they require multiple initial values to begin integration in each time slice and are not compatible with the proposed method in \cref{sec:GParareal}.} numerical integrators.
The first, referred to as the \emph{fine solver} $\F$, is a computationally expensive integrator that propagates an initial value at $t_j$, over an interval of length $\Delta T$, and returns a solution with high numerical accuracy at $t_{j+1}$.
In this paper, we assume that $\F$ provides sufficient numerical accuracy to solve \eqref{eq:ODE} for the solution to be considered `exact', i.e.\ $U_j = u(t_j)$.
The objective is to calculate the exact solutions
\begin{align}
    \label{eq:true_sol}
    U_{j} = \F(U_{j-1}) \quad \text{for} \quad j=1,\dots,J, \quad \text{where} \quad U_0 \defeq u_0,
\end{align}
\emph{without} running $\F$ $J$ times sequentially, as this calculation is assumed to be computationally intractable. To avoid this, parareal locates iteratively improved approximations $U^k_j$, where $k=0,1,2,\dots$ is the iteration number, that converge toward $U_j$ (note that $U^k_0=U_0=u_0 \ \forall k \geq 0$).
To do this, parareal uses a second numerical integrator $\G$, referred to as the \emph{coarse solver}.
$\G$ propagates an initial value at $t_j$ over an interval of length $\Delta T$, however, it has lower numerical accuracy and is computationally inexpensive to run compared to $\F$.
This means that $\G$ can be run serially across a number of time slices to provide relatively cheap low accuracy solutions whilst $\F$ is permitted only to run in parallel over multiple slices.

\subsection{The algorithm} \label{subsec:para_alg}
To begin (iteration $k=0$), approximate solutions to \eqref{eq:ODE} are calculated sequentially using $\G$, on a single processor, such that  
\begin{align} \label{eq:coarse}
    U^0_{j} = \G(U^0_{j-1}) \quad j=1,\dots,J.
\end{align}
Following this, the fine solver propagates each approximation in \eqref{eq:coarse} \emph{in parallel}, on $J$ processors, to obtain $\F(U^0_j)$ for $j=0,\dots,J-1$.
These values are then used (during iteration $k=1$) in the predictor-corrector  
\begin{align} \label{eq:pred-corr}
    U^{k}_{j} = \underbrace{\G(U^{k}_{j-1})}_{\text{predict}} + \underbrace{\F(U^{k-1}_{j-1}) - \G(U^{k-1}_{j-1})}_{\text{correct}} \quad \text{for} \quad j = 1,\dots,J.
\end{align}
Here, $\G$ is applied sequentially to predict the solution at the next time step, before being corrected by the residual between coarse and fine values found during the previous iteration (note that \eqref{eq:pred-corr} cannot be calculated in parallel).
This is a discretised approximation of the Newton--Raphson method for locating the true roots $U_j$ with initial guess \eqref{eq:coarse} \citep{gander2007}.
For a pre-defined tolerance $\varepsilon > 0$, the parareal solution $U^k_j$ is deemed to have converged up to time $t_I$ if 
\begin{align} \label{eq:stop}
	| U^k_j - U^{k-1}_j | < \varepsilon \quad \forall j \leq I.
\end{align}
This criterion is standard for parareal \citep{garrido2006,gander2008}, however, other criteria, e.g.\ taking the average relative error between fine solutions over a time slice \citep{samaddar2010,samaddar2019} or measuring the total energy of the system, could be used instead. 
Unconverged solution values, i.e.\ $U^k_j$ for $j > I$, are updated in future iterations ($k > 1)$ by initiating further parallel $\F$ runs on each $U^k_j$, followed by an update using \eqref{eq:pred-corr}.
The algorithm stops once $I=J$, converging in $k$ (out of $J$) iterations.
The version of parareal described here and implemented in \Cref{sec:numerics} does not iterate over solutions that have already converged, avoiding the waste of computational resources \citep{elwasif2011,pentland2022,garrido2006}.
Extending parareal to the full nonautonomous system in \eqref{eq:ODE} is straightforward:
see \cite{gander2007} for notation and \cite{pentland2022} for pseudocode.

\subsection{Convergence and computational complexity} \label{subsec:para_remarks}
After $k$ iterations, the solution states up to time $t_k$ (at minimum) have converged, as the exact initial condition ($u_0$) has been propagated by $\F$ at least $k$ times.
Therefore, if parareal converges in $k=J$ iterations, then the solution will be equal to the one found by calculating \eqref{eq:true_sol} serially, at an even higher computational cost. Convergence\footnote{For parareal to converge, the solvers $\F$ and $\G$ must satisfy specific mathematical conditions \citep{maday2005,bal2005}.} in $k \ll J$ iterations is necessary if significant parallel speed-up is to be realised.
Refer to \cite{gander2007,gander2008} for derivations of explicit parareal error bounds.

\begin{figure}[htbp]
    \centering
    \includegraphics[width=0.95\textwidth]{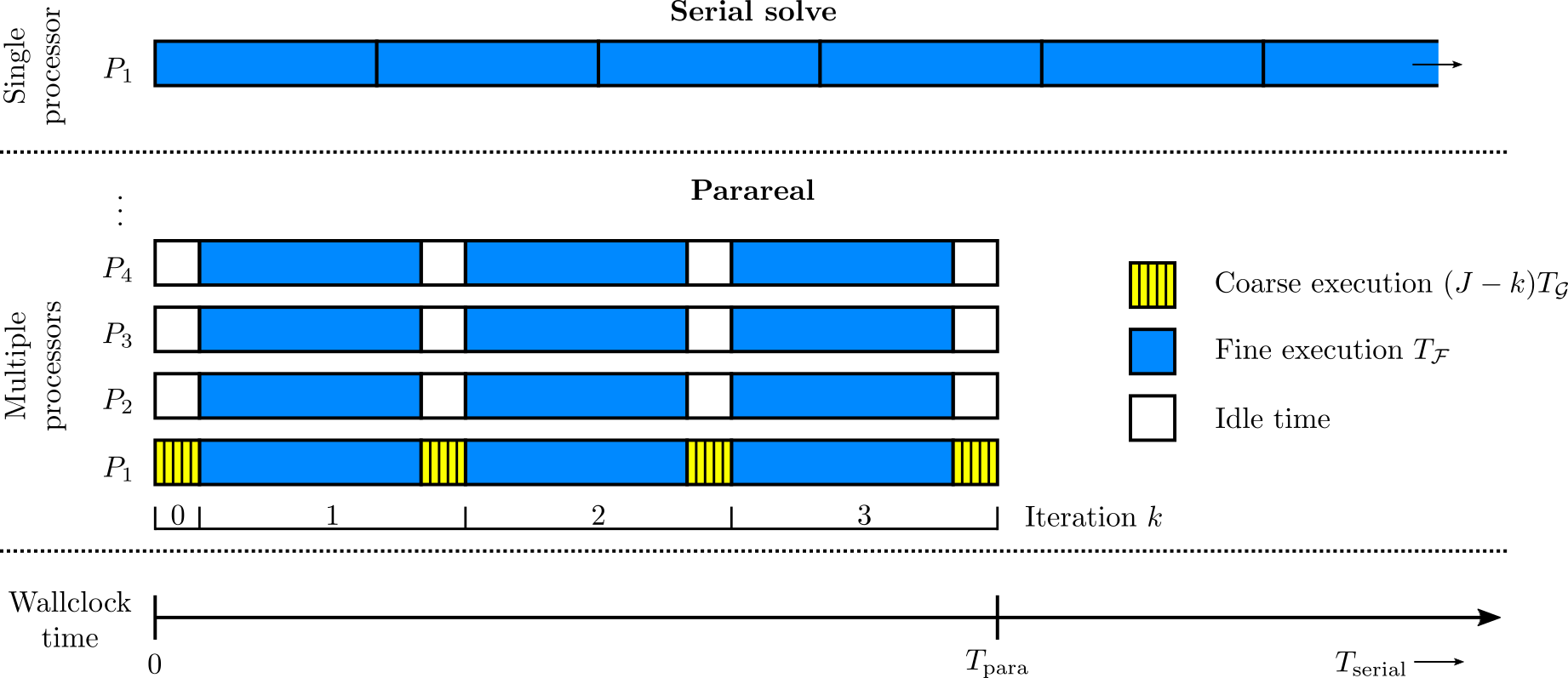}
    \caption{Computational task scheduling during three iterations of parareal compared with a full serial integration.
    The coloured blocks represent the wallclock time any given processor spends on a task.
    Coarse runs are shown in yellow, fine runs in blue, and any idle time in white.
    The wallclock time is given on the axis at the bottom, indicating both $T_{\text{para}}$ and $T_{\text{serial}}$.}
    \label{fig:processors}
\end{figure}

Without loss of generality, assume running $\F$ over any $[t_j,t_{j+1}]$, $j \in \{0,\dots,J-1 \}$, takes wallclock time $T_{\mathcal{F}}$ (denote time $T_{\mathcal{G}}$ similarly for $\G$). Therefore, calculating \eqref{eq:true_sol} using $\F$ serially, takes approximately $T_{\text{serial}} = J T_{\mathcal{F}}$ seconds. Using parareal, the total wallclock time (in the worst case, excluding any serial overheads) can be approximated by
\begin{align} \label{eq:Tpara}
    T_{\text{para}} \approx \underbrace{J T_{\mathcal{G}}}_{\text{Iteration 0}} + \sum_{i=1}^{k} \underbrace{\bigl( T_{\mathcal{F}} + (J-i) T_{\mathcal{G}} \bigr)}_{\text{Iterations 1 to $k$}} =  k T_{\mathcal{F}} + (k+1) \left( J - \frac{k}{2} \right) T_{\mathcal{G}}.
\end{align}
The approximate parallel speed-up is therefore
\begin{align} \label{eq:Spara}
    S_{\text{para}} \approx \frac{T_{\text{serial}}}{T_{\text{para}}} = \left[ \frac{k}{J} + (k+1) \left( 1-\frac{k}{2J} \right) \frac{T_{\mathcal{G}}}{T_{\mathcal{F}}} \right]^{-1}.
\end{align}
To maximise \eqref{eq:Spara}, both the convergence rate $k$ and the ratio $T_{\mathcal{G}}/T_{\mathcal{F}}$ should be as small as possible.
In practice, however, there is a trade-off between these two quantities as fast $\G$ solvers (with sufficient accuracy to still guarantee convergence) typically require more iterations to converge, increasing $k$.
An illustration of the computational task scheduling during the first few iterations of parareal vs.\ a full serial integration is given in \cref{fig:processors}---optimised scheduling of parareal is studied in \cite{elwasif2011}.

Selecting a fast but accurate coarse solver remains a trial and error process, entirely dependent on the system being solved. Typically, $\G$ is chosen such that it has a coarser temporal resolution/lower numerical accuracy \citep{samaddar2010,farhat2003,baffico2002,trindade2006}, a coarser spatial resolution (when solving PDEs) \citep{samaddar2019,ruprecht2014}, and/or uses simplified model equations \citep{engblom2009,legoll2020,meng2020} compared to $\F$.
In \cref{sec:GParareal}, we aim to widen the pool of choices for $\G$ by using a GP emulator to capture variability in the residual $\F-\G$ and showcase its effectiveness by demonstrating that GParareal can converge to a solution in cases where parareal cannot in \cref{sec:numerics}.

\section{GParareal} \label{sec:GParareal}
In this section, we present the GParareal algorithm, in which a GP emulator is used in the analogue of parareal's predictor-corrector step.
Suppose we seek the same high resolution numerical solutions to \eqref{eq:ODE} as expressed in \eqref{eq:true_sol}, denoted now as $V_j$ instead of $U_j$. Furthermore, we denote the iteratively improved approximations from GParareal as $V^k_j$ (as before, $V^k_0 = V_0 = u_0 \ \forall k \geq 0$).

In parareal, the predictor-corrector \eqref{eq:pred-corr} updates the numerical solutions at iteration $k$ using a correction term based on information calculated during the \emph{previous} iteration $k-1$. 
We propose the following update rule, again based on a coarse prediction and multi-fidelity correction, that instead refines solutions using information from the \emph{current} iteration $k$, rather than $k-1$:
\begin{align} \label{eq:trick}
    V^k_{j} = \F(V^k_{j-1}) &= (\F - \G + \G)(V^k_{j-1}) \nonumber \\
    &= \underbrace{(\F - \G) (V^k_{j-1})}_{\text{correction}} + \underbrace{\G(V^k_{j-1})}_{\text{prediction}} \quad 1 \leq k < j \leq J.
\end{align}
If $V^k_{j-1}$ is known, the prediction is rapidly calculable, however the correction is not known explicitly without running $\F$ at expensive cost. 
We propose using a GP emulator to model this correction term, trained on \emph{all} previously obtained evaluations of $\F$ and $\G$.
The emulator returns a Gaussian distribution over $(\F - \G) (V^k_{j-1})$ from which we can extract an explicit value and carry out the refinement in \eqref{eq:trick}. 

In \Cref{subsec:new_alg}, we present the algorithm, giving an explanation of the kernel hyperparameter optimisation process in \Cref{subsec:hyperparams} and providing error analysis in \Cref{subsec:convergence}.
In \Cref{subsec:complexity}, we detail the computational complexity, remarking that given large enough runtimes for the fine solver, an iteration of GParareal runs in approximately the same wallclock time as parareal.
Again, to simplify notation, we first detail GParareal for an autonomous scalar-valued ODE, i.e.\ $\bm{f}(t,\bm{u}(t)) \defeq f(u(t))$ in \eqref{eq:ODE}.
The extension to the multivariate nonautonomous case is described in \Cref{subsec:vectorised}.

\newpage
\subsection{The algorithm} \label{subsec:new_alg}

\subsubsection*{Gaussian process emulator}
Before solving \eqref{eq:ODE}, we define a GP prior \citep{rasmussen2006} over the unknown correction function $\F-\G$.
This function maps an initial value $x_j \in \mathcal{U}$ at time $t_j$ to the residual difference between $\F(x_j)$ and $\G(x_j)$ at time $t_{j+1}$.
More formally, we define the GP prior
\begin{equation} \label{eq:GP_prior}
 \F-\G \sim \mathcal{GP} ( m , \kappa ),
\end{equation}
with mean function $m \colon \mathcal{U} \to \Reals$ and covariance kernel $\kappa\colon \mathcal{U} \times \mathcal{U} \to \Reals$.
Given some vectors of initial values, $\bm{x},\bm{x}' \in \mathcal{U}^J$, the corresponding vector of means is denoted $\mu (\bm{x}) = ( m(x_j) )_{j=0,\dots,J-1}$ and the covariance matrix $K(\bm{x},\bm{x}') = ( \kappa(x_i,x'_j) )_{i,j=0,\dots,J-1}$.
The correction term is expected to be small, depending on the accuracy of both $\F$ and $\G$, hence we define a zero-mean process, i.e.\ $m(x_j) = 0$.
Ideally, the covariance kernel will be chosen based on any prior knowledge of the solution to \eqref{eq:ODE}, e.g.\ regularity/smoothness.
If no information is available \emph{a priori} to simulation, we are free to select any appropriate kernel.
In this work, we use the square exponential (SE) kernel 
\begin{align} \label{eq:kernel}
    \kappa(x,x') = \sigma^2 \exp \left( -\frac{(x-x')^2}{2\ell^2} \right) , \quad \text{for some} \quad x,x' \in \mathcal{U}.
\end{align}
The kernel hyperparameters, denoting the output length scale $\sigma^2$ and input length scale $\ell^2$, are referred to collectively in the vector $\bm{\theta}$ and need to be initialised prior to simulation.
The algorithm proceeds as follows;
see \Cref{algorithm} for pseudocode.

\subsubsection*{Iteration $k=0$}
Firstly, run $\G$ sequentially from the exact initial value, on a single processor, to locate the coarse solutions
\begin{equation} \label{eq:G_solve}
    V^0_{j} = \G(V^0_{j-1}) \quad j = 1,\dots,J.
\end{equation}
Store these solutions in the vector $\bm{x} \defeq (V^0_0,\dots,V^0_{J-1})^\intercal$ for use in the GP emulator. 

\subsubsection*{Iteration $k=1$}
Use $\F$ to propagate the values in \eqref{eq:G_solve} on each time slice in \emph{parallel}, on $J$ processors, to obtain the following values at $t_{j}$ 
\begin{equation} \label{eq:F_solve}
    \F(V^0_{j-1}) \quad j = 1,\dots,J.
\end{equation}
At this stage, we diverge from the parareal method. 
Given $\bm{x}$, store the values of $\F-\G$, using \eqref{eq:G_solve} and \eqref{eq:F_solve}, in the vector
\begin{equation} \label{eq:data}
    \bm{y} \defeq \bigl( (\F - \G)(x_j) \bigr)_{j=0,\dots,J-1}.
\end{equation}
At this point, the inputs $\bm{x}$ and evaluations $\bm{y}$ are used to optimise the kernel hyperparameters $\bm{\theta}$ via maximum likelihood estimation---see \Cref{subsec:hyperparams}. 
Conditioning the prior \eqref{eq:GP_prior} using the acquisition data $\bm{x}$ and $\bm{y}$, the GP posterior over $(\F-\G)(x')$, where $x' \in \mathcal{U}$ is some initial value in the state space, is given by
\begin{equation} \label{eq:GP_post}
(\F-\G)(x') \ | \ \bm{x},\bm{y} \sim \mathcal{N} \bigl( \hat{\mu}(x'), \hat{K}(x',x') \bigr),
\end{equation}
with mean
\begin{equation} \label{eq:GP_post_mean}
    \hat{\mu}(x') = \underbrace{ \mu(x')}_{=0} + K(x',\bm{x}) [K(\bm{x},\bm{x})]^{-1} \bigl( \bm{y} - \underbrace{\mu(\bm{x})}_{=\bm{0}} \bigr)
\end{equation}
and variance
\begin{equation} \label{eq:GP_post_cov}
    \hat{K}(x',x') = K(x',x') - K(x',\bm{x})[K(\bm{x},\bm{x})]^{-1} K(\bm{x},x').
\end{equation}

Now we wish to determine updated solutions $V^1_j$ at each mesh point.
Given $\F$ has been run once, the exact solution is known at time $t_1$.
Specifically, at $t_0$ we know $V^k_0 = V_0 \ \forall k \geq 0$ and at $t_1$ we know $V^k_1 = V_1 = \F(V^1_0) \ \forall k \geq 1$.
At $t_2$, the exact solution $V_2 = \F(V^1_1)$ is unknown, hence we need to calculate its value without running $\F$ again.
To do this, we re-write the exact solution using \eqref{eq:trick}:
\begin{align}
    \label{eq:v2}
    V^1_2 = \F(V^1_1) = (\F - \G + \G)(V^1_1) = \underbrace{(\F - \G) (V^1_1)}_{\text{correction}} + \underbrace{\G(V^1_1)}_{\text{prediction}}.
\end{align}
Both terms in \eqref{eq:v2} are initially unknown, but the prediction can be calculated rapidly at low computational cost while the correction can be inferred using the GP posterior \eqref{eq:GP_post} with $x'=V^1_1$.
Therefore, we obtain a Gaussian distribution over the solution
\begin{align} \label{eq:V2}
    V^1_2 \sim \mathcal{N} \bigl( \hat{\mu}(V^1_1) + \G(V^1_1), \hat{K}(V^1_1,V^1_1) \bigr),
\end{align}
with variance stemming from uncertainty in the GP emulator.
Repeating this process to determine a distribution for the solution at $t_3$ by attempting to propagate the random variable $V^1_2$ using $\G$ is computationally infeasible for nonlinear IVPs.
To tackle this and be able to propagate $V_2^1$, we approximate the distribution by taking its mean value,
\begin{align*}
    V^1_2 = \hat{\mu}(V^1_1) + \G(V^1_1).
\end{align*}
This approximation is a convenient way of minimising computational cost, at the price of ignoring uncertainty in the GP emulator---see \Cref{sec:discussion} for a discussion of possible alternatives.

\begin{figure}[tbp]
    \centering
    \includegraphics[width=0.99\textwidth]{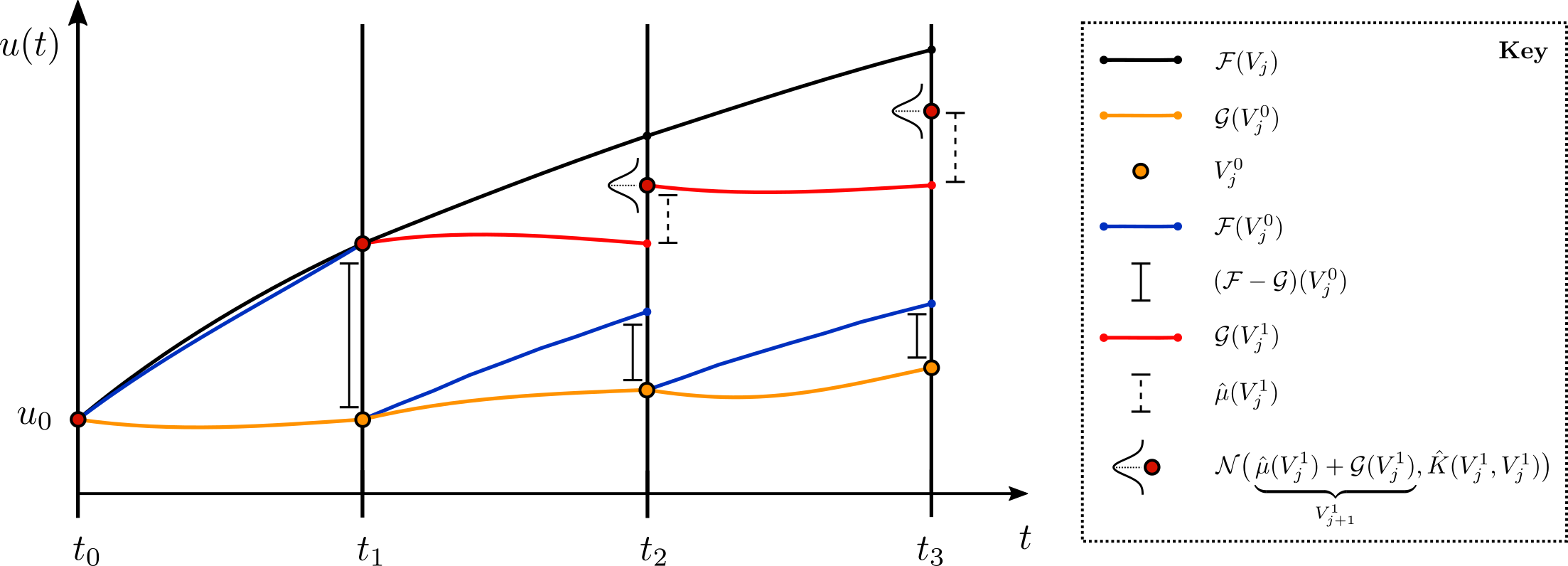}
    \caption{Schematic of the first iteration of GParareal.
    The `exact' solution over $[t_0,t_3]$ is shown in black, with the first coarse and fine (parallel) runs given in yellow and blue respectively.
    Solid bars represent the residual between these solutions \eqref{eq:data}.
    The predictions, i.e.\ the second coarse runs, are shown in red and the corresponding corrections from the GP emulator are represented by the dashed bars.
    The updated solutions \eqref{eq:Vj_expec} at the end of the iteration are represented by the red dots.
    Note the black and blue lines in $[t_0.t_1]$ should overlap but have been made not to for clarity.}
    \label{fig:GParareal}
\end{figure}

The update process, applying \eqref{eq:trick} and then approximating the Gaussian distribution by taking its expectation, is repeated sequentially for later $t_j$, yielding the approximate solutions
\begin{align} \label{eq:Vj_expec}
    V^1_{j} = \hat{\mu}(V^1_{j-1}) + \G(V^1_{j-1}) \quad \text{for} \quad j = 3,\dots,J.
\end{align}
This process is illustrated in \Cref{fig:GParareal}.
Finally, we impose stopping criteria \eqref{eq:stop}, identifying which $V^1_j$ for $j\leq I$ have converged.
Using the same stopping criteria as parareal will allow us to compare the performance of both algorithms in \Cref{sec:numerics}. 

\subsubsection*{Iteration $k \geq 2$}
If the stopping criteria is not met, i.e.\ $I < J$, we can iteratively update any unconverged solutions by re-applying the previous steps.
This means calculating $\F(V^{k-1}_j)$, $j = I,\dots,J-1$, in parallel and then storing new evaluations $\hat{\bm{y}} = \bigl( (\F-\G)(V^{k-1}_j) \bigr)^{\intercal}_{j=I,\dots,J-1}$, with corresponding inputs $\hat{\bm{x}} = (V^{k-1}_I,\dots,V^{k-1}_{J-1})^{\intercal}$. Hyperparameters are then re-optimised and the GP is re-conditioned using \emph{all} prior acquisition data, i.e.\ $\bm{x} = [\bm{x};\hat{\bm{x}}]$ and $ \bm{y} = [\bm{y};\hat{\bm{y}}]$, generating an updated posterior.
Here, $[\bm{a};\bm{b}]$ denotes the vertical concatenation of column vectors $\bm{a}$ and $\bm{b}$.
The update rule is then applied such that we obtain
\begin{align*}
    V^k_{j} = \hat{\mu}(V^k_{j-1}) + \G(V^k_{j-1}) \quad \text{for} \quad j = I+2,\dots,J.
\end{align*}
Once $I=J$, the solution, the number of iterations $k$ taken to converge, and the acquisition data $\bm{x}$ and $\bm{y}$ are returned.
A key advantage of GParareal is that the acquisition data can be used in future GParareal simulations (as ``legacy data'') to provide the GP emulator with more data and therefore exploit additional speedup---this will be demonstrated in \cref{sec:numerics}. 

\subsection{Kernel hyperparameter optimisation} \label{subsec:hyperparams}
The hyperparameters $\bm{\theta}$ of the kernel $\kappa$ will need to be optimised in light of the acquisition data $\bm{y}$ (and corresponding input data $\bm{x}$). 
We optimise each element of $\bm{\theta}$ such that it maximises its (log) marginal likelihood \citep{rasmussen2004}. 
To do this, first define $g(x) \defeq (\F-\G)(x)$ and $\bm{g} \defeq (g(x_j))_{j=0,\dots,N-1}$, where $N$ is the current length of $\bm{x}$ (and $\bm{y}$). 
Given the evaluations $\bm{y}$ are noise-free, the likelihood of obtaining such data is $p(\bm{y}|\bm{g},\bm{x},\bm{\theta}) = \delta(\bm{y} - \bm{g})$, where $\delta(\cdot)$ is the multidimensional Dirac delta function. 
The marginal likelihood, given $\bm{x}$ and $\bm{\theta}$, is therefore
\begin{align*}
    p(\bm{y}|\bm{x},\bm{\theta}) &= \int \underbrace{p(\bm{y}|\bm{g},\bm{x},\bm{\theta})}_{\text{likelihood}} \underbrace{p(\bm{g}|\bm{x},\bm{\theta})}_{\text{prior}} \, \rd \bm{g} 
    = \int \delta(\bm{y} - \bm{g}) \mathcal{N}(\bm{g}|\bm{0},K(\bm{x},\bm{x})) \, \rd \bm{g}
    =  \mathcal{N}(\bm{y}|\bm{0},K(\bm{x},\bm{x})),
\end{align*}
where $\mathcal{N}(\bm{y}|\bm{0},K(\bm{x},\bm{x}))$ denotes the probability density function of a multivariate Gaussian distribution evaluated at $\bm{y}$, with mean vector $\bm{0}$ and covariance matrix $K(\bm{x},\bm{x})$ that depends on $\bm{\theta}$, see \eqref{eq:kernel}.
The hyperparameters in $\bm{\theta}$ can then be estimated numerically by maximising the log marginal likelihood using any gradient-based optimiser.
Optimisation is carried out once per iteration (up until the hyperparameters do not change significantly between iterations) and hyperparameters from the prior iteration are used as to start the optimisation at the current iteration. 


\subsection{Error Analysis} \label{subsec:convergence}
In this section, we are interested in analysing the absolute error
\begin{equation}\label{abserr}
    e^k_j \defeq | V_j - V^k_j |,
\end{equation}
between the exact and GParareal solution at iteration $k$ and time $t_j$.
We show that this error has an upper bound proportional to the \emph{fill distance} (defined below) of the dataset at iteration $k$.
To do this, we now denote the input dataset at iteration $k$ as $\bm{x}_k$ rather than $\bm{x}$ (because the dataset size strictly increases with each iteration of GParareal) and, similarly, denote the output dataset $\bm{y}$ as $\bm{y}_k$. 
We also introduce some assumptions on the solvers $\F$ and $\G$, and a known result on the consistency of the GP posterior mean $\hat{\mu}$ \eqref{eq:GP_post_mean} to the true correction function $g = \F - \G$.
For clarity, we re-state the GParareal update rule
\begin{align} \label{eq:update}
    V^{k}_{j} = \G(V^{k}_{j-1}) + \hat{\mu}(V^{k}_{j-1}), \quad 1 \leq k < j \leq J.
\end{align}

\subsubsection{Preparatory assumptions and results}
First, we state some assumptions on $\F$ and $\G$, as in \cite{gander2008}.
\begin{assumption}[Exact fine solver $\F$] \label{assump:exactflow}
$\F$ solves \eqref{eq:ODE} exactly such that
\begin{align} \label{eq:exact}
    V_j = \F (V_{j-1}), \quad j = 1,\ldots,J.
\end{align}
\end{assumption}

\begin{assumption}[One-step coarse solver $\G$] \label{assump:coarseapprox}
$\G$ is a one-step numerical solver with uniform local truncation error $\mathcal{O}(\Delta T^{p+1})$, for $p \geq 1$, such that
\begin{align*}
    \F(u) - \G(u) = c_1 (u) \Delta T^{p+1} + c_2 (u) \Delta T^{p+2} + \ldots,
\end{align*}
for $u \in \Reals$ and continuously differentiable functions $c_i(u)$, $i = 1,2,\ldots$. 
For $u, v \in \Reals$, we can then write
\begin{align} \label{eq:FG2}
    | \left(\F(u) - \G(u)\right) - \left(\F(v) - \G(v)\right) | \leq C_1 \Delta T^{p+1} | u - v |,
\end{align}
where $C_1>0$ is the Lipschitz constant for $c_1(u)$.
\end{assumption}

\begin{assumption}[Lipschitz coarse solver $\G$] \label{assump:lipschitz}
$\G$ satisfies the Lipschitz condition
\begin{align} \label{eq:lipschitz}
     | \G(u) - \G(v) | \leq L_{\mathcal{G}} | u - v |,
\end{align}
for $u, v \in \Reals$ and some $L_{\mathcal{G}} > 0$.
\end{assumption}

Next, we define the concepts required to state a result on the consistency of the GP posterior mean.
Firstly, we define the \emph{fill distance} $h_{\bm{x}_k}$ to be the largest smallest distance between any point $v \in \mathcal{U}$ and any point $v_i \in \bm{x}_k$, i.e.\ 
\begin{align*}
    h_{\bm{x}_k} \defeq \sup_{v \in \mathcal{U}} \inf_{v_i \in \bm{x}_k} | v - v_i |.
\end{align*}
It should be clear that each data point $v_i \in \bm{x}_k$ is also contained in $\mathcal{U}$.
Secondly, we define the \emph{reproducing kernel Hilbert space} (RKHS), a Hilbert space $H_{\kappa}(\mathcal{U})$ of functions $g\colon \mathcal{U} \rightarrow \Reals$ with inner product $\langle \cdot,\cdot \rangle_{H_{\kappa}(\mathcal{U})}$.
See \citet{stuart2018} for a more formal definition and conditions on the inner product itself. 
We can now state the following result on the GP posterior mean consistency, adapted from \citet[Theorem~11.14]{wendland2004}.

\begin{theorem}[GP posterior mean consistency] \label{thm:GP_consist}
Suppose $\mathcal{U} \subset \Reals$ is a bounded interval and let $\kappa$ be the SE kernel.
Denote the GP posterior mean, built using $\bm{x}_k$, $\bm{y}_k$, and $\kappa$ \eqref{eq:GP_post_mean} as $\hat{\mu}$ and the function being emulated as $g \in H_{\kappa}(\mathcal{U})$.
Then, for every $\tau \in \mathbb{N}$, there exist constants $h_0(\tau)$ and $C_{\tau} > 0$ such that
\begin{align*}
    | g(v) - \hat{\mu}(v) | \leq C_{\tau} h^{\tau}_{\bm{x}_k} | g |_{H_{\kappa}(\mathcal{U})} \quad \forall v \in \mathcal{U},
\end{align*}
provided that $h_{\bm{x}_k} \leq h_0(\tau)$.
Note that  $| g |_{H_{\kappa}(\mathcal{U})}^2 = \langle g,g \rangle_{H_{\kappa}(\mathcal{U})}$.
\end{theorem}
See \citet[Theorem~11.14]{wendland2004} for a more general version of this result that holds when $\mathcal{U} \subset \Reals^d$ and for derivatives of both $g$ and $\hat{\mu}$.
It should be noted that \cref{thm:GP_consist} holds when $g \in H_{\kappa}(\mathcal{U})$, i.e.\ the function of interest lies within the RKHS of the SE kernel. 
If this is not the case, convergence issues may arise (see \cite{karvonen2022a,karvonen2022b}) and one would need to choose an alternative kernel function.
For consistency results involving Mat{\'e}rn kernels, see \cite{stuart2018}.

\subsubsection{Error bound for GParareal solutions}

\begin{theorem}[GParareal error bound] \label{thm:GParareal_error}
Suppose the solvers used in GParareal satisfy \Cref{assump:exactflow,assump:coarseapprox,assump:lipschitz}, and that the conditions required for \cref{thm:GP_consist} hold.
Then, the absolute error \eqref{abserr} of the GParareal solution to the autonomous scalar-valued ODE, i.e.\ $\bm{f}(t,\bm{u}(t)) \defeq f(u(t))$ in \eqref{eq:ODE}, at iteration $k$ and time $t_j$ satisfies
\begin{align*}
e^k_j \leq 
\begin{dcases}
        \Lambda_k \sum_{i=0}^{j-(k+1)} A^i & \ 1 \leq k < j \leq J, \\[0.5em]
        0 & \ 0 \leq j \leq k \leq J.
\end{dcases}
\end{align*}
where $A = C_1 \Delta T^{p+1} + L_{\mathcal{G}}$ and $\Lambda_k = C_{\tau} h^{\tau}_{\bm{x}_k} |g|_{H_{\kappa}(\mathcal{U})}$.
\end{theorem}
\begin{proof}
First, consider the case $0 \leq j \leq k \leq J$.
For $j=0$, recall that $V^k_0 = V_0 \ \forall k \geq 0$ by definition, hence $e^k_0 = 0 \ \forall k \geq 0$.
For $j=1$, we seek $V^1_1 = \F(V^1_0)$ which we in fact know from applying $\F$ to $V^0_0$ during the prior iteration (i.e. $k=0$). 
Therefore, we have that
\begin{align*}
    V^1_1 = \F(V^1_0) = \F(V^0_0) = \F(V_0) = V_1 \quad &\Rightarrow \quad V^k_1 = V_1 \ \forall k \geq 1 \quad \Rightarrow \quad e^k_1 = 0 \ \forall k \geq 1.
\end{align*}
We can repeat this process iteratively up to $j=J$ to show that
\begin{align*}
    V^J_J = \F(V^J_{J-1}) = \F(V^{J-1}_{J-1}) = \F(V_{J-1}) = V_J \quad &\Rightarrow \quad V^k_J = V_J \ \forall k \geq J \\
    &\Rightarrow \quad e^k_J = 0 \ \forall k \geq J.
\end{align*}
Now, consider the case $1 \leq k < j \leq J$. 
Using the update rule \eqref{eq:update}, that $\F$ is the exact solver \eqref{eq:exact}, and adding and subtracting the terms $g(V^k_j)$ and $\G(V_j)$, we can write 
\begin{align*}
    e^k_{j+1} &= | V_{j+1} - V^k_{j+1} | = | \F(V_j) - \big( \G(V^k_j) + \hat{\mu}(V^k_j) \big) | \\
    &= | \F(V_j) - \big( \G(V^k_j) + \hat{\mu}(V^k_j) \big) \pm g(V^k_j) \pm \G(V_j)|.
\end{align*}
Applying the triangle inequality and the definition of $g$, we obtain
\begin{multline*}
    e^k_{j+1} \leq | \big( \F(V_j) - \G(V_j) \big) - \big( \F(V^k_j) - \G(V^k_j) \big) | \\
    + | \G(V_j) - \G(V^k_j) | + | g(V^k_j) - \hat{\mu}(V^k_j) |.
\end{multline*}
On the right hand side, the first term can be bounded using \eqref{eq:FG2}, the second by \eqref{eq:lipschitz}, and the third using \cref{thm:GP_consist}, yielding the recursion
\begin{align*}
    e^k_{j+1} &\leq A e^k_j  + \Lambda_k,
\end{align*}
where $A = C_1 \Delta T^{p+1} + L_{\mathcal{G}}$ and $\Lambda_k = C_{\tau} h^{\tau}_{\bm{x}_k} |g|_{H_{\kappa}(\mathcal{U})}$.
This recursion can be solved using the initial condition $e^k_j = 0 \ \forall k \geq j$ to obtain the desired result.
\end{proof}

\Cref{thm:GParareal_error} shows that the error is proportional to the fill distance at iteration $k$ and that GParareal will recover the exact solution at time $t_j$ after $k=j$ iterations.

Note that this result is rather general in the sense that we consider the fill distance with respect to the entire space $\mathcal{U} \subset \Reals$, whereas in reality we would measure the fill distance with respect to a moderately sized compact interval $\mathcal{V} \subset \mathcal{U}$ in which the solution $u(t)$ lies $\forall t \in [t_0,T]$.
Essentially, the accuracy of the GP posterior mean outside of $\mathcal{V}$ is inconsequential to the GParareal scheme because the mean will never be evaluated outside of $\mathcal{V}$.
Also note, the result will generalise for GParareal applied to systems of ODEs by using norms and the generalised version of \cref{thm:GP_consist} in \cite{wendland2004}.

\subsection{Computational complexity} \label{subsec:complexity}
The complexity of GParareal can be calculated similarly to that of parareal---refer back to \cref{subsec:para_remarks} for notation.
In GParareal, an additional cost is incurred when (serially) conditioning the emulator on acquisition/legacy data and optimising the hyperparameters.
During the $k$th iteration, up to $kJ$ evaluations of $\F-\G$ have been collected, hence standard cubic complexity GP conditioning scales like $\mathcal{O}(k^3 J^3)$ in terms of floating point operations (and $\mathcal{O}(k^2 J^2)$ per hyperparameter).
Given a fixed number of time slices $J$, let $T_{\text{GP}}(k)$ represent the total wallclock time taken to condition and optimise hyperparameters of the GP (using up to $kJ$ observations) at iteration $k$---note this is a strictly increasing function of $k$.
Ignoring serial overheads, we can write down the total wallclock time for GParareal as
\begin{align} 
    T_{\text{GPara}} \approx J T_{\mathcal{G}} + \sum_{i=1}^{k} \bigl( T_{\mathcal{F}} + (J-i) T_{\mathcal{G}} + T_{\text{GP}}(i) \bigr) \nonumber\\
    = k T_{\mathcal{F}} + (k+1) \left( J - \frac{k}{2} \right) T_{\mathcal{G}} + T_{\text{GP}}, \label{eq:TGPara} 
\end{align}
where $T_{\text{GP}} := \sum_{i=1}^k T_{\text{GP}}(i)$.
The approximate parallel speed-up is then given by 
\begin{align} \label{eq:SGPara}
    S_{\text{GPara}} \approx \left[ \frac{k}{J} + (k+1) \left( 1-\frac{k}{2J} \right) \frac{T_{\mathcal{G}}}{T_{\mathcal{F}}} + \frac{1}{J} \frac{T_{\text{GP}}}{T_{\mathcal{F}}} \right]^{-1}.
\end{align}
Therefore, in addition to the parareal requirements that $k$ be as small as possible and $T_{\mathcal{G}} \ll T_{\mathcal{F}}$, GParareal requires that $T_{\text{GP}} \ll T_{\mathcal{F}}$ in order to maximise parallel speedup.
If this is the case, the complexity of GParareal is approximately the same as parareal.

This suggests that if $k$ and/or $J$ are large, then the cost of the emulation may begin to dominate that of the fine solver, limiting the parallel speedup from GParareal, see \cref{sec:numerics}.
This, however, need not hinder the usability of GParareal for a number of reasons.
Firstly, time-parallelisation is typically deployed on problems where additional parallel speedup is needed beyond that achieved by traditional domain decomposition, i.e.\ spatio-temporal PDEs.
This means that if $P$ processors are required for the space-parallel computations of the PDE and $J$ processors for the time-parallel computations, then $JP$ processors are required in total.
For moderate to large values of $P$, only leftover HPC resources are available to exploit time-parallelism and so $J$ typically cannot be chosen very large, somewhat limiting how large $T_{\text{GP}}$ can be.
Secondly, in the scenario that both $T_{\text{GP}}$ and $T_{\mathcal{F}}$ are small, one does not need to use a time-parallel method in the first place, as $\F$ can simply be run serially in this case. 
Thirdly, if both $T_{\text{GP}}$ and $T_{\mathcal{F}}$ are large or of a similar order, then one can reduce $T_{\text{GP}}$ by reducing the number of time slices $J$, thereby increasing $T_{\mathcal{F}}$ at the same time. 

Whilst there is no way to control the final value of $k$ obtained by either parareal or GParareal, there are ways of reducing $T_{\text{GP}}$ using more efficient non-cubic complexity, emulation methods. For example, one could make use of sparse GPs, parallel matrix inversion methods, or sparse approximate linear algebra techniques \citep{schaefer2021} to reduce the cost of evaluating the inverse kernel matrix $K(\bm{x},\bm{x})^{-1}$.
One could also reduce $T_{\text{GP}}$ by clustering the input data points and training `local' GPs in parallel \citep{snelson2007} or instead use inducing points to average over input data points that are located close together in state space \citep{quinonero2005,snelson2006}---see \cite{murphy2023} for additional methods.
To reduce the, often significant, cost of hyperparameter optimisation, one may deploy parallel optimisation routines if available or, as we implement in \cref{sec:numerics}, stop the optimisation once additional data no longer improves the hyperparameter estimates.

\subsection{Generalisation to ODE systems} \label{subsec:vectorised}
The methodology in \Cref{subsec:new_alg} can be generalised to solve systems of $d$ autonomous ODEs.
Accordingly, the correction term we wish to emulate is now vector-valued, i.e.\ $\mathcal{U} \subset \Reals^d$, hence we require a vector-valued (or multi-output) GP, rather than a scalar GP. 

The simplest approach is to model each output of $\F-\G$ independently, whereby we use $d$ scalar GPs (sharing the same vector-valued inputs in state space) to emulate each output. 
This requires initialising $d$ GP emulators, each with their own covariance kernel $\kappa_i$ (usually the same for consistency) and corresponding hyperparameters $\bm{\theta}_i$---to be optimised independently using their own respective observation datasets $\bm{y}^{(i)}$, $i=1,\dots,d$.
In this case, the $d$ GP emulators can be conditioned/optimised independently of one another and so we make use of the idle processors to carry out these computations in an embarrassingly parallel fashion to reduce the total GP complexity from $\mathcal{O}(d k^3 J^3)$ to $\mathcal{O}(k^3 J^3)$ each iteration---the same as the scalar case.

The more complex approach is to jointly emulate the outputs of $\F-\G$ by modelling cross-covariances between outputs via the method of co-kriging \citep{cressie1993}. 
A number of co-kriging techniques exist (see \cite{alvarez2011} for a brief overview), one of which is the linear model of coregionalisation that models the joint, block-diagonal, covariance prior using a linear combination of the separate kernels $\kappa_i$.
Prior testing revealed that using this method did not improve performance enough to justify the added complexity, $\mathcal{O}(d^3 k^3 J^3)$ vs.\ $\mathcal{O}(d k^3 J^3)$ in the independent setting (results not reported). 
Some applications may require correlated output dimensions, hence we note the methodology here for any interested readers. 

As a final note, to solve nonautonomous systems of equations, i.e.\ \eqref{eq:ODE}, there are two possible approaches. 
One way is to include the time variable as an extra input to each of the $d$ scalar GPs---this requires a more carefully selected covariance kernel.
The other way is to re-write the $d$-dimensional nonautonomous system as a system of $d+1$ autonomous equations and solve as described above---this is the method we use in \Cref{sec:numerics}.


\section{Numerical experiments} \label{sec:numerics}
In this section, we present numerical experiments to compare the performance of GParareal and parareal on a number of low-dimensional ODE systems, namely the FitzHugh--Nagumo model, the chaotic R{\"ossler} system, a nonautonomous system, and the double pendulum system.
MATLAB code for GParareal, parareal, and the GP emulator as used in the experiments of this section can be found at \url{https://github.com/kpentland/GParareal}.

For simplicity, $\F$ and $\G$ are chosen to be explicit Runge--Kutta methods (RK) of order $q,p \in \{ 1,2,4,8 \}$, respectively ($q \geq p$).
Let $N_{\mathcal{F}}$ and $N_{\mathcal{G}}$ denote the number of time steps each solver uses over $[t_0,T]$.
For these experiments we built our own cubic complexity GP emulator to highlight the effectiveness of GParareal using standard out-the-box methods, postponing the implementation of more efficient and sophisticated emulation methods to a future work.
In the multivariate setting (recall \cref{subsec:vectorised}), we use a scalar output GP emulator (with isotropic SE covariance kernel) to model each output dimension of $\F-\G$ and assign each one its own processor, reducing the GP emulation costs by a factor of $d$.
Hyperparameter optimisation is carried out at each iteration, stopping when the (maximal) absolute difference between hyperparameters is larger than $10^{-2}$. 
The experiments are run on up 512 CPUs.

\subsection{FitzHugh--Nagumo model} \label{subsec:FHN}
\begin{figure}[b!]
    \centering
    \begin{subfigure}{0.47\linewidth}
        \includegraphics[width=\textwidth]{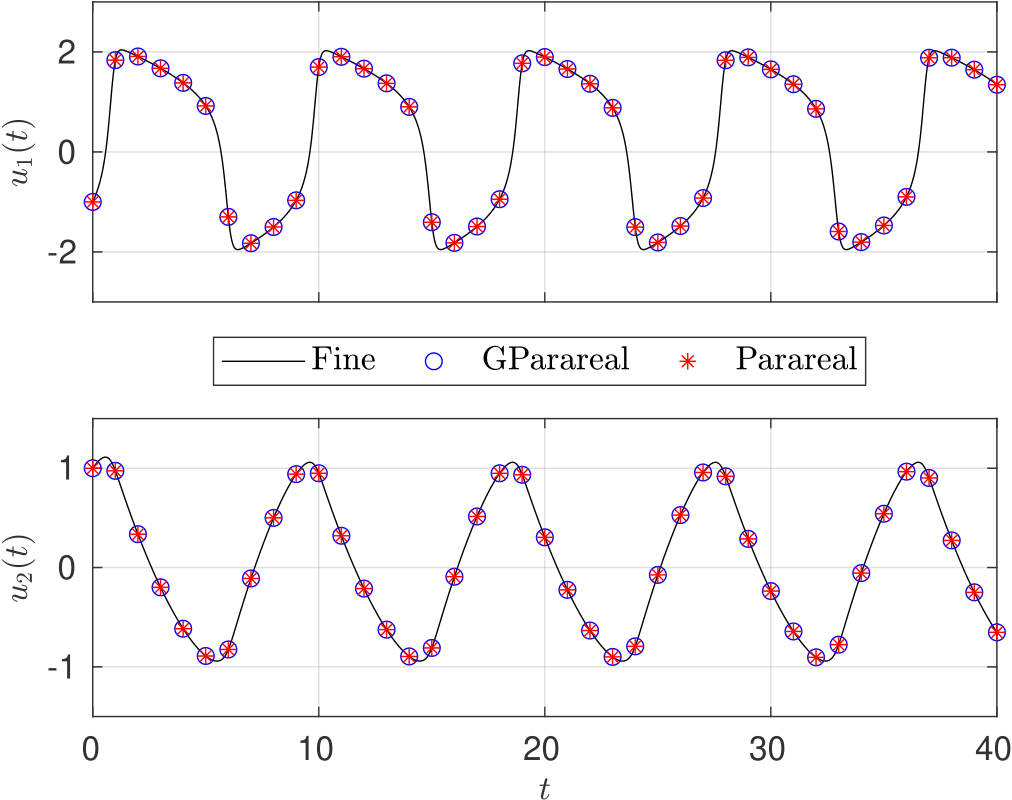}
        \caption{~}
    \end{subfigure}
    \begin{subfigure}{0.49\linewidth}
        \includegraphics[width=\textwidth]{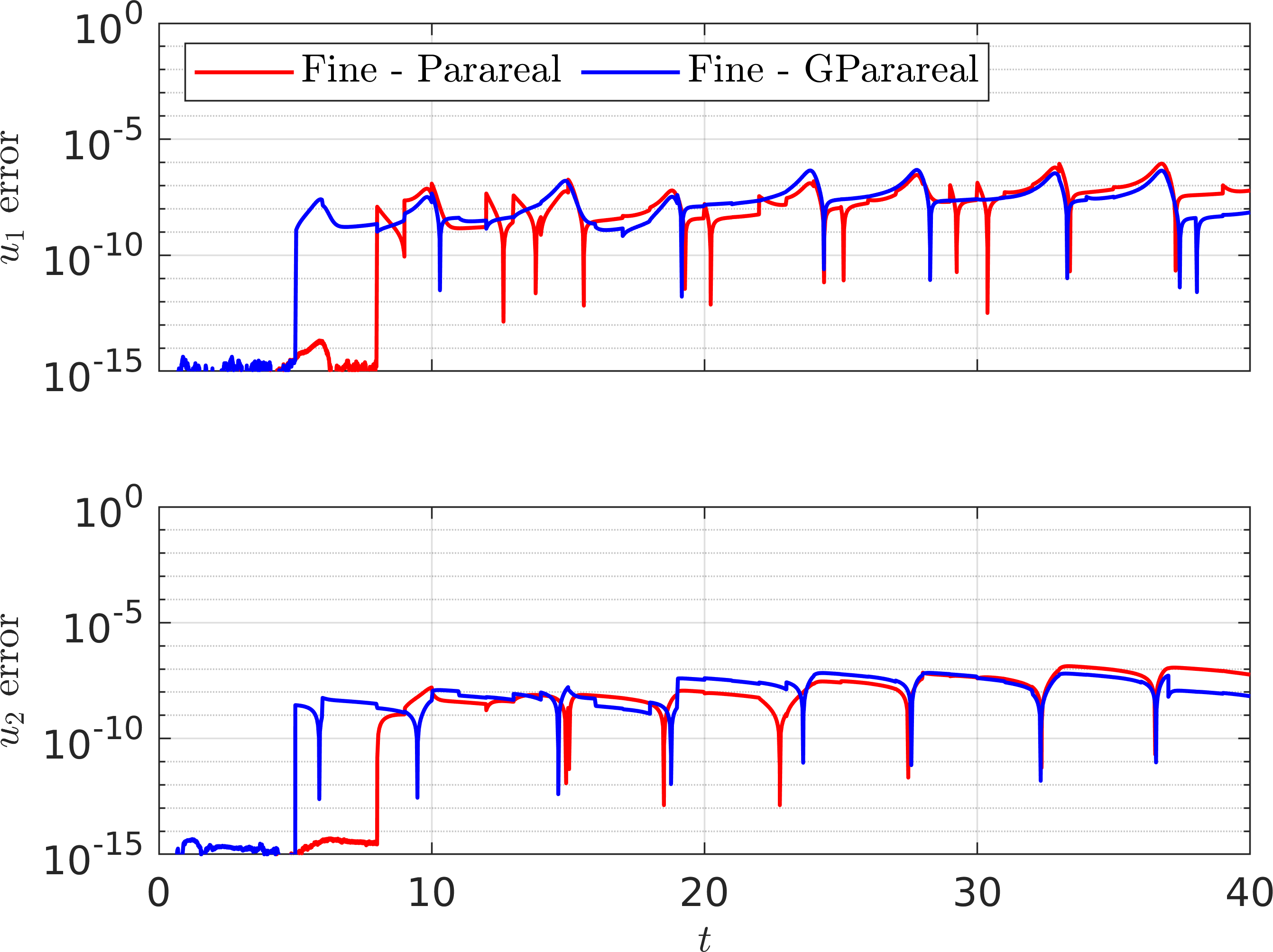}
        \caption{~}
    \end{subfigure}
        \begin{subfigure}{0.49\linewidth}
        \includegraphics[width=\textwidth]{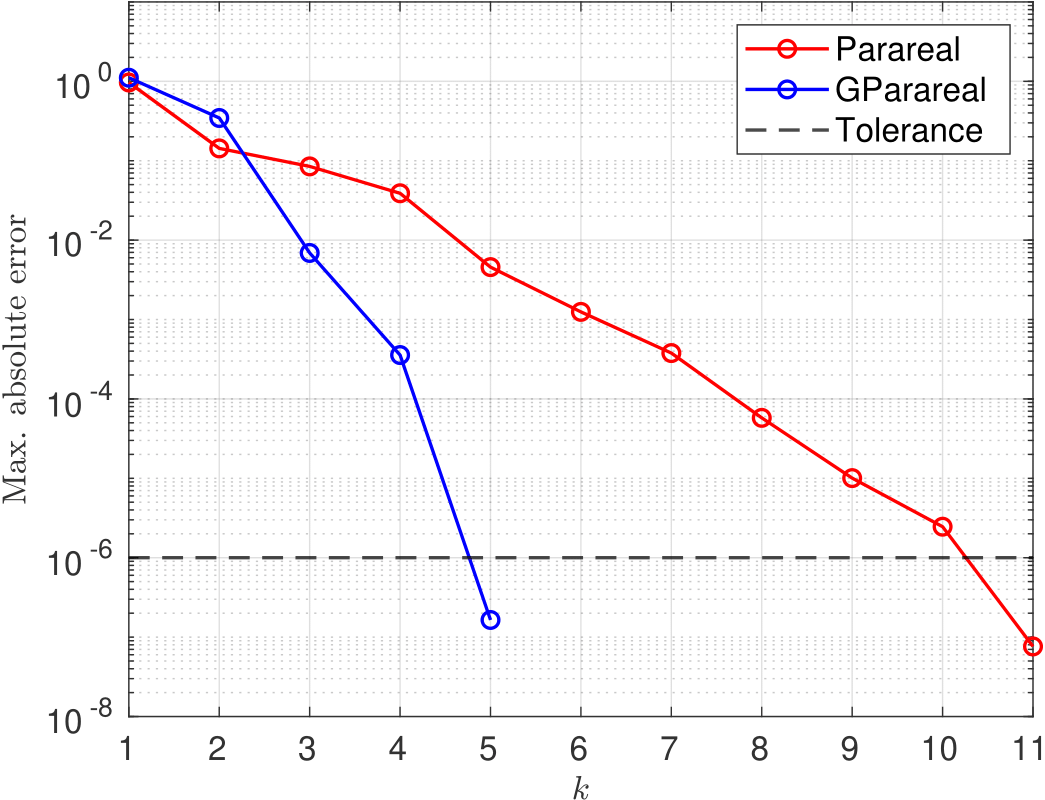}
    \caption{~}
    \end{subfigure}
        \begin{subfigure}{0.485\linewidth}
        \includegraphics[width=\textwidth]{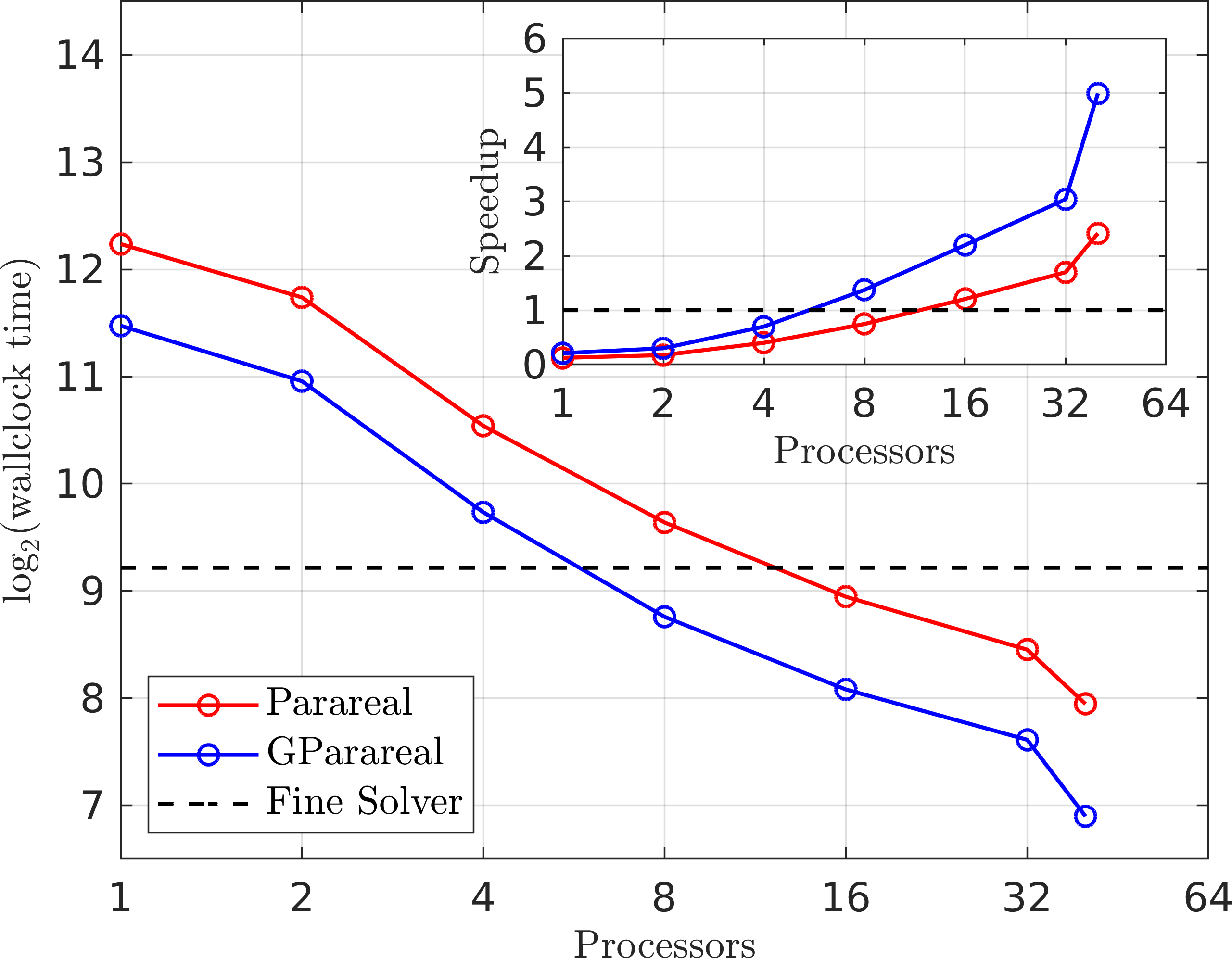}
    \caption{~}
    \end{subfigure}
    \caption{Numerical results obtained solving the FHN model \eqref{eq:FHN} for $\bm{u}_0 = (-1,1)^\intercal$. 
    (a) Time-dependent solutions using the fine solver, GParareal, and parareal---both GParareal and parareal plotted only at times $\bm{t}$ for clarity. 
    (b) The corresponding absolute errors between solutions from GParareal and parareal vs.\ the fine solution. 
    (c) Maximum absolute errors \eqref{eq:stop} of each algorithm at successive iterations $k$ until tolerance $\varepsilon = 10^{-6}$ is met. 
    (d) Median wallclock times (taken over 5 runs) of both algorithms against the number of processors (up to 40). The inset plot displays the corresponding parallel speedup.}
    \label{fig:FHN}
\end{figure}
In this experiment, we consider the FitzHugh--Nagumo (FHN) model \citep{fitzhugh1961,nagumo1962a} given by
\begin{equation} \label{eq:FHN}
    \frac{\rd u_1}{\rd t} = c \bigl(u_1 - \frac{u_1^3}{3} + u_2 \bigr), \quad \frac{\rd u_2}{\rd t} = -\frac{1}{c}(u_1 - a + bu_2),
\end{equation}
where we fix parameters $(a,b,c) = (0.2,0.2,3)$.
We integrate \eqref{eq:FHN} over $t \in [0,40]$, dividing the interval into $J=40$ slices, and set the tolerance for both GParareal and parareal to $\varepsilon = 10^{-6}$.
We use solvers $\G=\text{RK2}$ and $\F=\text{RK4}$ with $N_{\mathcal{G}} = 160$ and $N_{\mathcal{F}} = 1.6 \times 10^{5}$ steps respectively.

\begin{figure}[t!]
    \centering
    \begin{subfigure}{0.49\linewidth}
        \includegraphics[width=\textwidth]{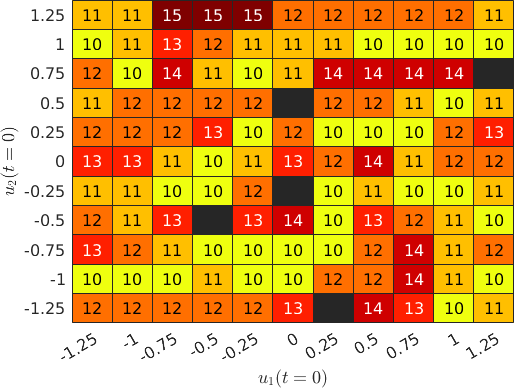}
        \caption{Parareal}
    \end{subfigure}
    \begin{subfigure}{0.49\linewidth}
        \includegraphics[width=\textwidth]{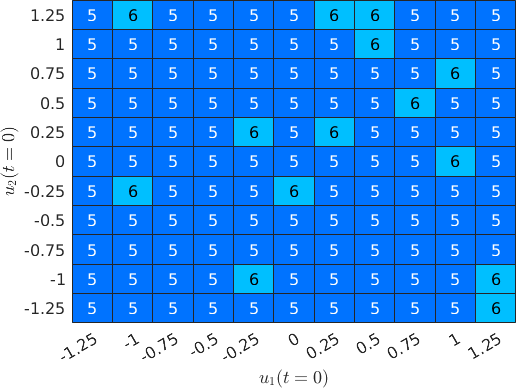}
        \caption{GParareal}
    \end{subfigure}
    \caption{Heat maps displaying the number of iterations taken until convergence $k$ of (a) parareal and (b) GParareal when solving the FHN model \eqref{eq:FHN} for different initial values $\bm{u}_0 \in [-1.25,1.25]^2$. Black boxes indicate where parareal returned a \texttt{NaN} value during simulation.}
    \label{fig:FHN_heat}
\end{figure}
\begin{figure}[b!]
    \centering
    \begin{subfigure}{0.475\linewidth}
        \includegraphics[width=\textwidth]{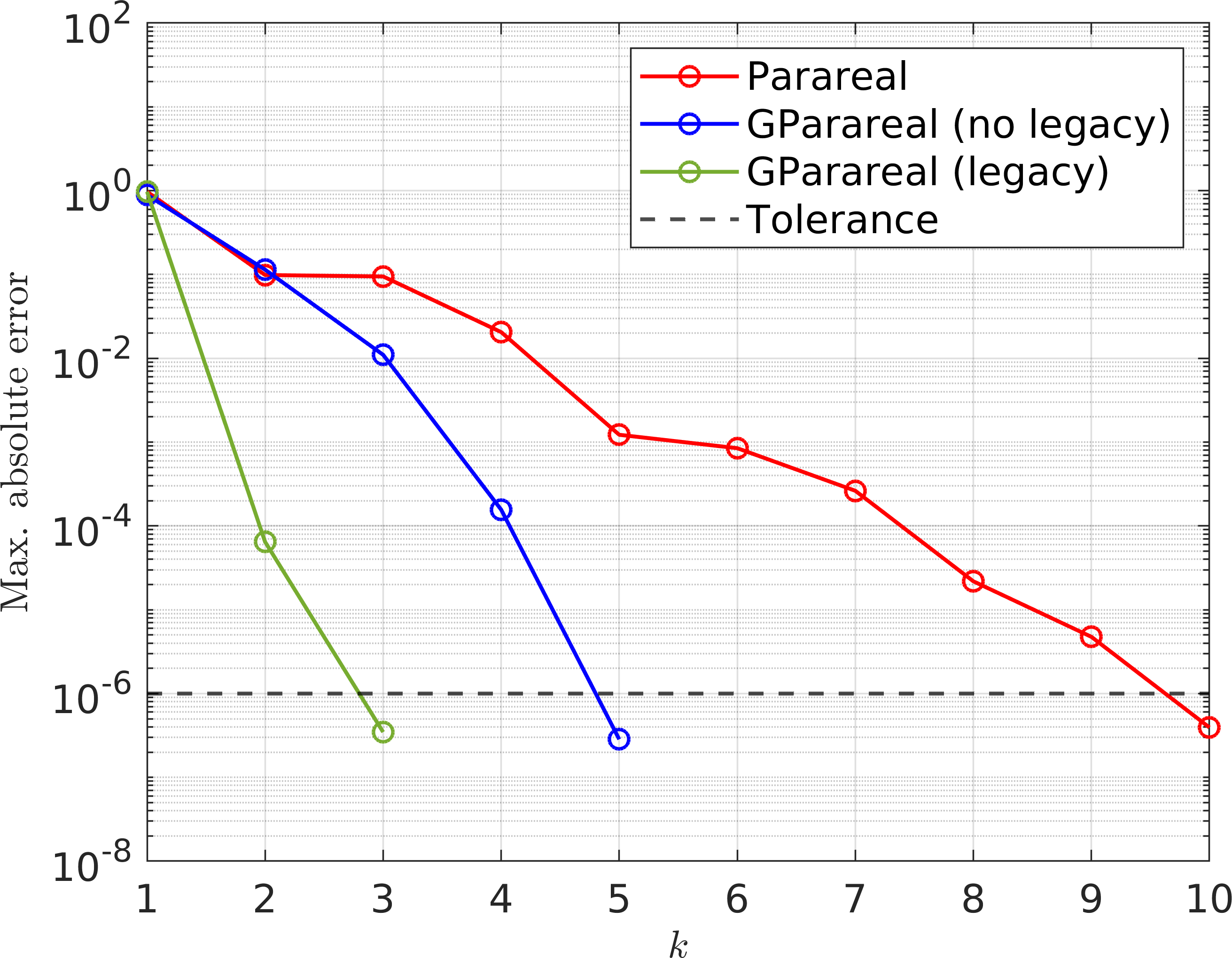}
        \caption{~}
    \end{subfigure}
    \begin{subfigure}{0.49\linewidth}
        \includegraphics[width=\textwidth]{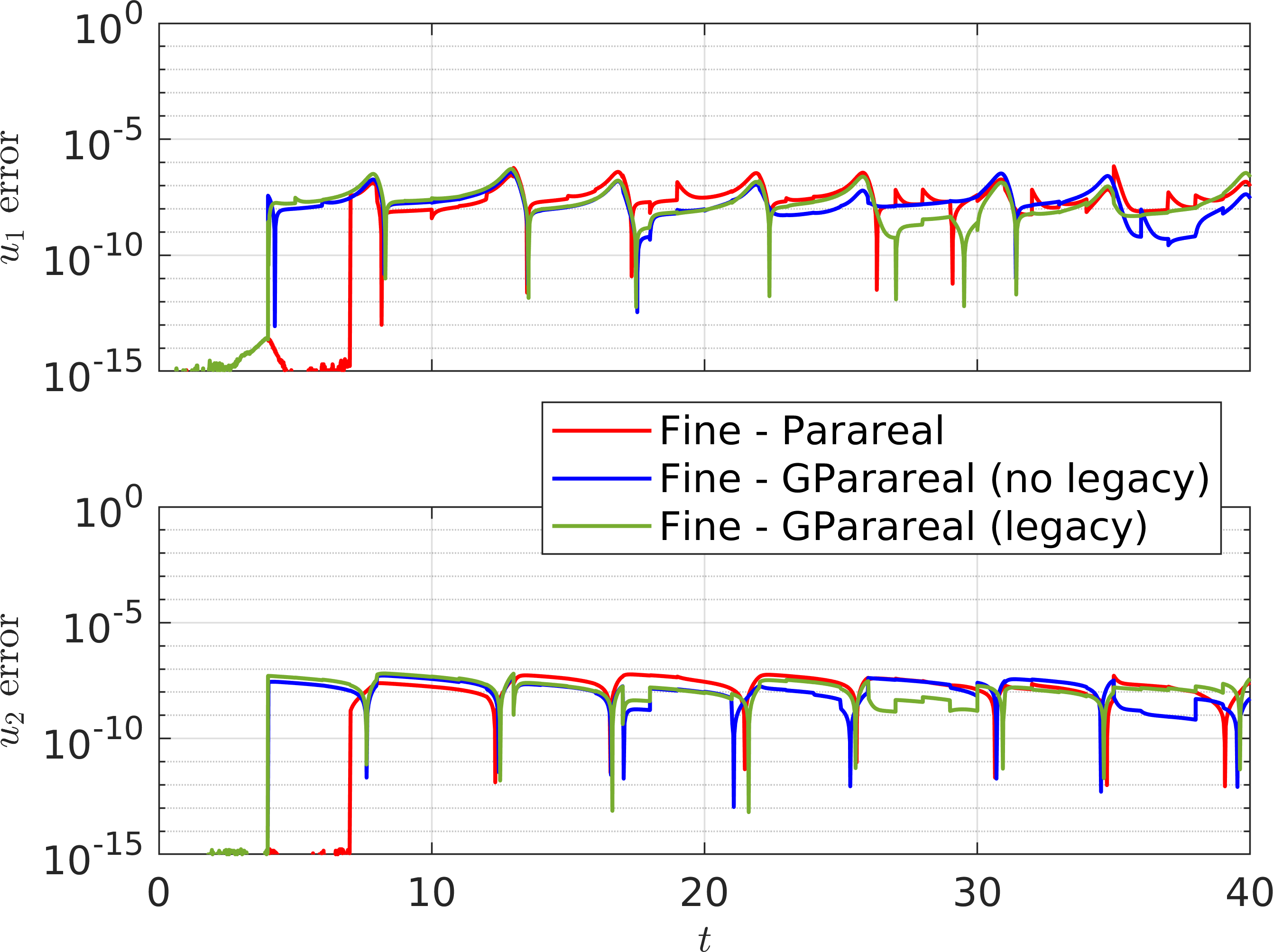}
        \caption{~}
    \end{subfigure}
    \caption{Numerical simulations solving \eqref{eq:FHN} for $\bm{u}(0) = (0.75,0.25)^\intercal$ using GParareal with and without access to legacy data, i.e.\ $\F-\G$ data obtained solving \eqref{eq:FHN} for  $\bm{u}(0) = (-1,1)^\intercal$. The parareal simulation of the same problem is also shown for comparison. (a) Maximum absolute errors \eqref{eq:stop} against iteration number $k$ until tolerance $\varepsilon=10^{-6}$ met. (b) Time-dependent errors of the corresponding numerical solutions from each simulation vs.\ the fine solution.}
    \label{fig:FHN_legacy}
\end{figure}

In \cref{fig:FHN}(a), we solve \eqref{eq:FHN} with initial condition $\bm{u}_0 = (-1,1)^\intercal$ using both algorithms.
Observe that the accuracy of GParareal is of approximately the same order as the solution obtained using parareal---when comparing both to the serially obtained fine solution (\cref{fig:FHN}(b)).
Note, however, that in \cref{fig:FHN}(c), GParareal takes six fewer iterations to converge to these solutions than parareal does.
As a result, GParareal locates a solution in faster wallclock time than parareal, see \cref{fig:FHN}(d), with an almost $5$-fold speedup vs.\ the serial solver---over twice the $2.4$-fold speedup obtained by parareal.
Note that we increase $N_{\mathcal{F}}$ to $1.6 \times 10^{8}$ to ensure $\F$ is expensive to run and realise parallel speedup in \cref{eq:FHN}(d) (as both algorithms require $T_{\mathcal{G}}/T_{\mathcal{F}} \ll 1$).

To compare the convergence of both methods more broadly, we solve \eqref{eq:FHN} for a range of initial values.
The heatmap in \cref{fig:FHN_heat}(a) illustrates how the convergence of parareal is highly dependent, not just on the solvers in use, but also the initial values at $t=0$, taking anywhere from 10 to 15 iterations to converge.
For some initial values, parareal does not converge at all, with solutions blowing up (returning \texttt{NaN} values) due to the low accuracy of $\G$.
In direct contrast, see \cref{fig:FHN_heat}(b), GParareal converges more quickly and more uniformly due to the flexibility provided by the emulator, taking just five or six iterations to reach tolerance for all the initial values tested. 
This demonstrates how using an emulator can enable convergence even when $\G$ has poor accuracy.

\begin{figure}[t!]
    \centering
    \includegraphics[width=0.49\textwidth]{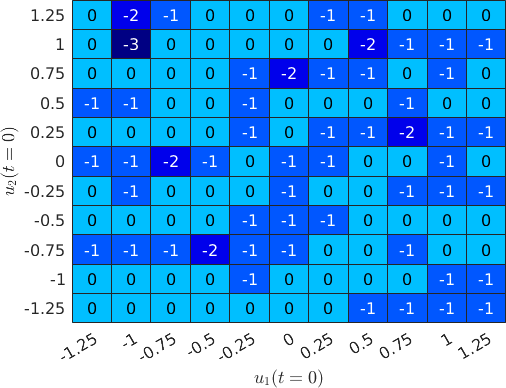}
    \caption{Heat map displaying the decrease in the number of iterations taken until convergence of GParareal when solving \eqref{eq:FHN} for different initial values $\bm{u}_0 \in [-1.25,1.25]^2$ \emph{with} legacy data compared to without, i.e.\ compared to \cref{fig:FHN_heat}(b). Legacy data was obtained by solving \eqref{eq:FHN} for $\bm{u}_0 = (-1,1)^\intercal$.}
    \label{fig:FHN_heat_legacy}
\end{figure}

Until now, GParareal simulations have been carried out using only acquisition data.
In \cref{fig:FHN_legacy}, we demonstrate how GParareal can use both acquisition and legacy data to converge in fewer iterations than without the legacy data.
Approximately $kJ = 5 \times 40 = 200$ legacy data points, obtained solving \eqref{eq:FHN} for $\bm{u}_0 = (-1,1)^\intercal$, are stored and made available to the GP emulator when solving \eqref{eq:FHN} for alternate initial values $\bm{u}_0 = (0.75,0.25)^\intercal$.
In \cref{fig:FHN_legacy}(a), we can see that convergence takes two fewer iterations with the legacy data than without.
Accuracy of the solutions obtained from these simulations is again shown to be of the order of the parareal solution in both cases---see \cref{fig:FHN_legacy}(b).
Repeating the experiment from \Cref{fig:FHN_heat}(b) with the same legacy data for a range of initial values we see that $k$ is either unchanged or improved in all cases, see \cref{fig:FHN_heat_legacy}.
It should be noted that conditioning the GP and optimising hyperparameters using the legacy data comes at extra (serial) computational cost and checks should be made to ensure that $T_{\mathcal{F}} \gg T_{\text{GP}}$. 
These results illustrate that using GParareal (with or without legacy data) we can solve and evaluate the dynamics of the FHN model in significantly lower wallclock time than parareal. 

\subsection{R{\"o}ssler system} \label{subsec:rossler}
\begin{figure}[b!]
    \centering
    \begin{subfigure}{0.46\linewidth}
        \includegraphics[width=\textwidth]{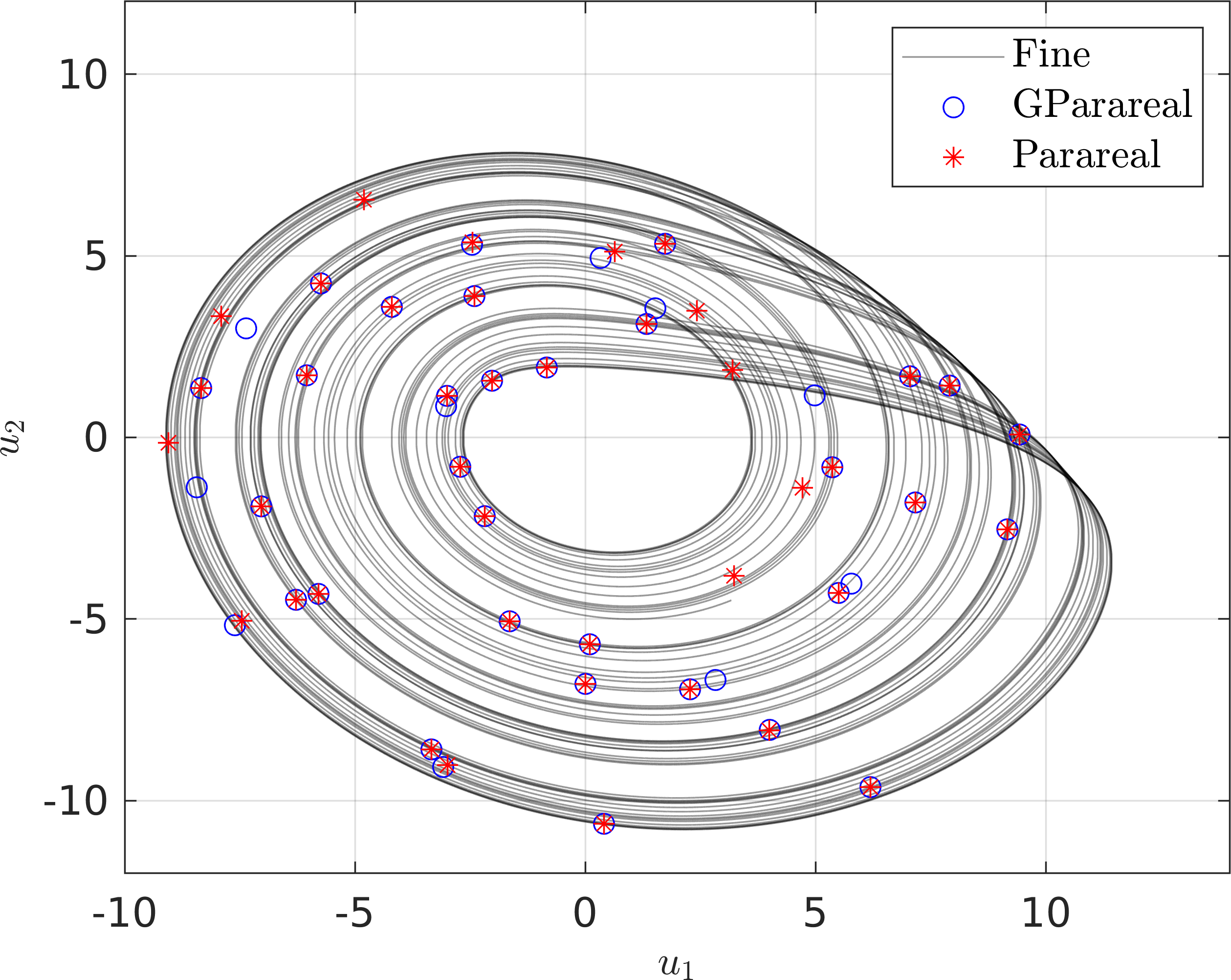}
        \caption{~}
    \end{subfigure}
    \begin{subfigure}{0.49\linewidth}
        \includegraphics[width=\textwidth]{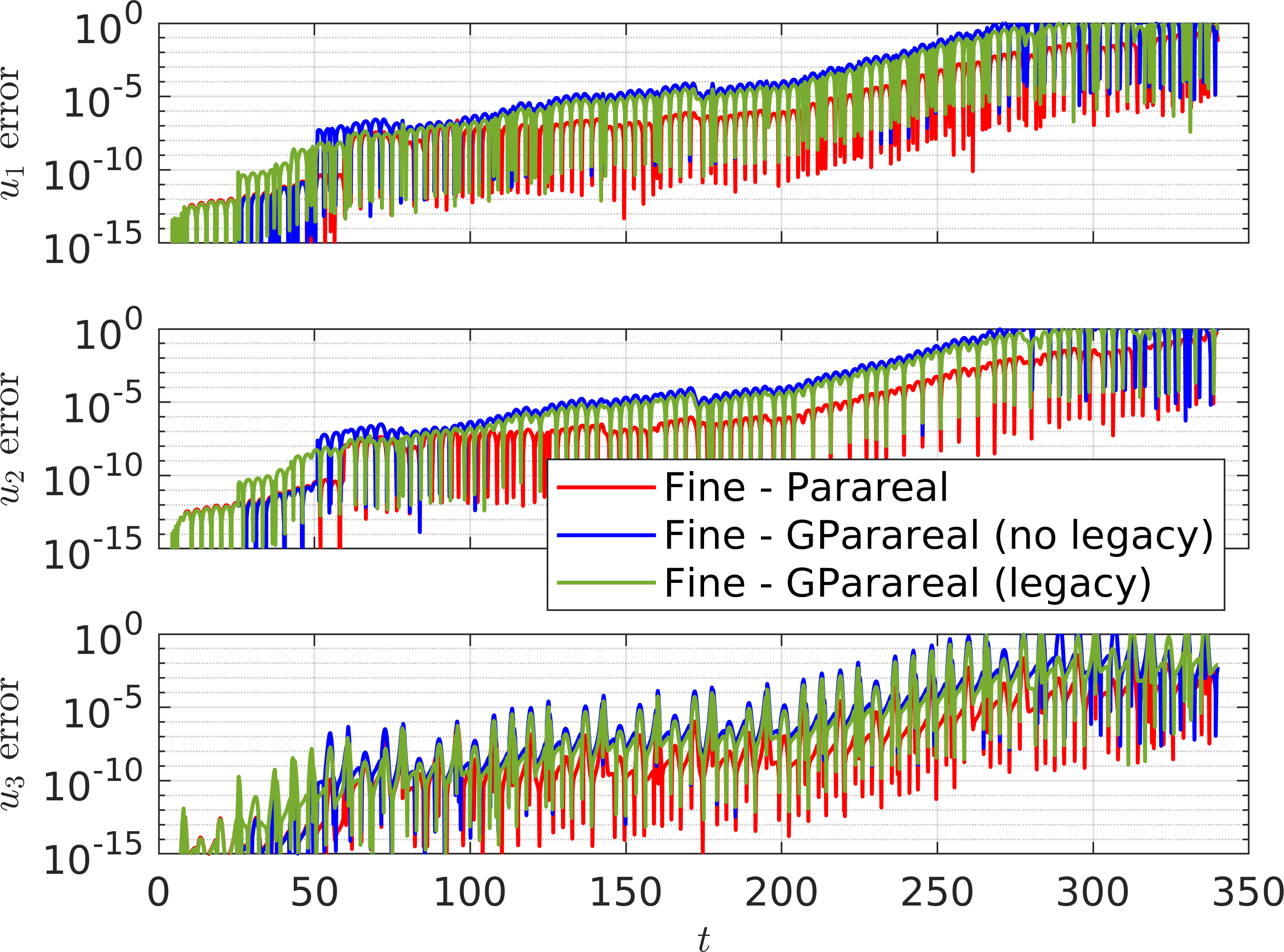}
        \caption{~}
    \end{subfigure}
    \begin{subfigure}{0.49\linewidth}
        \includegraphics[width=\textwidth]{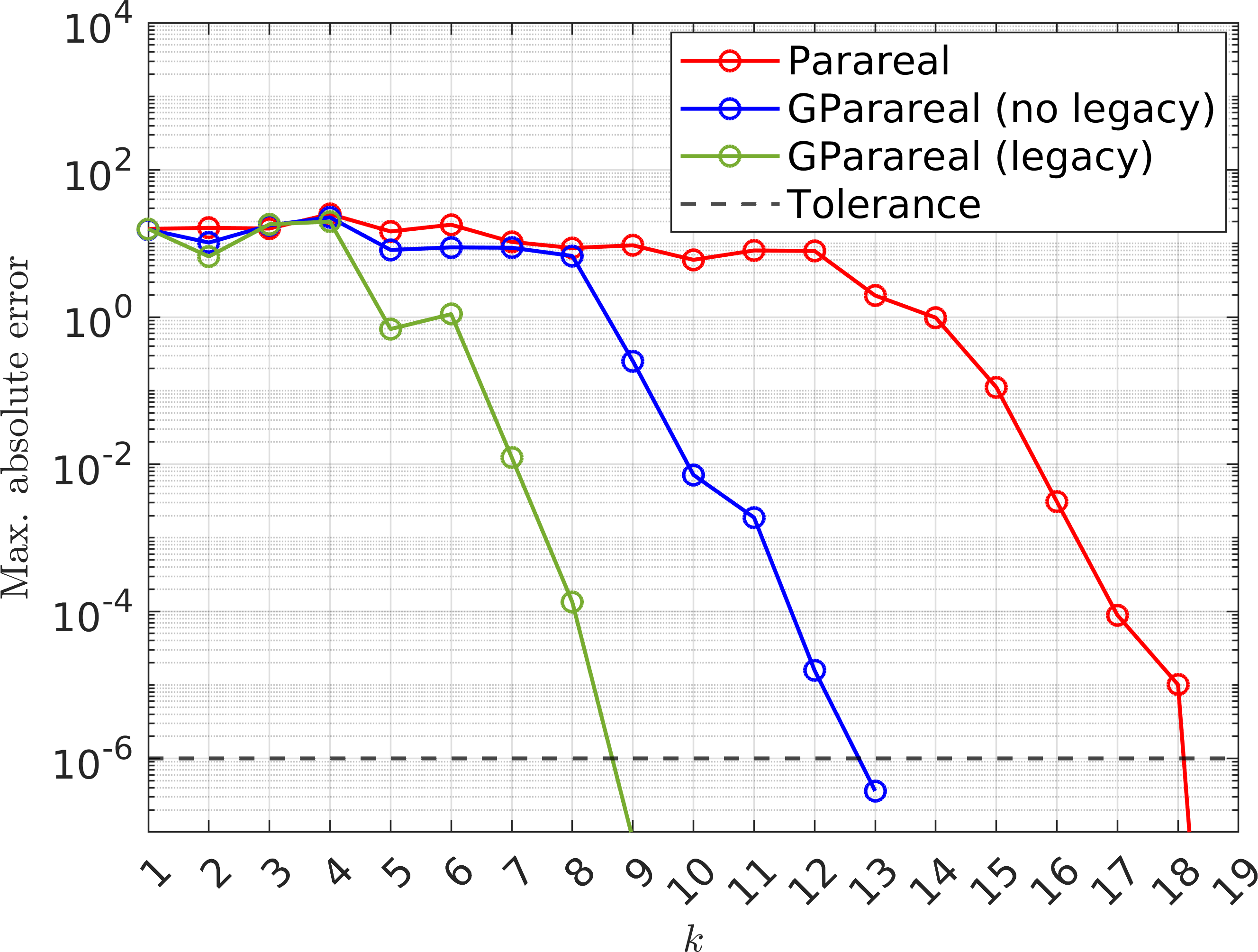}
        \caption{~}
    \end{subfigure}
    \begin{subfigure}{0.475\linewidth}
        \includegraphics[width=\textwidth]{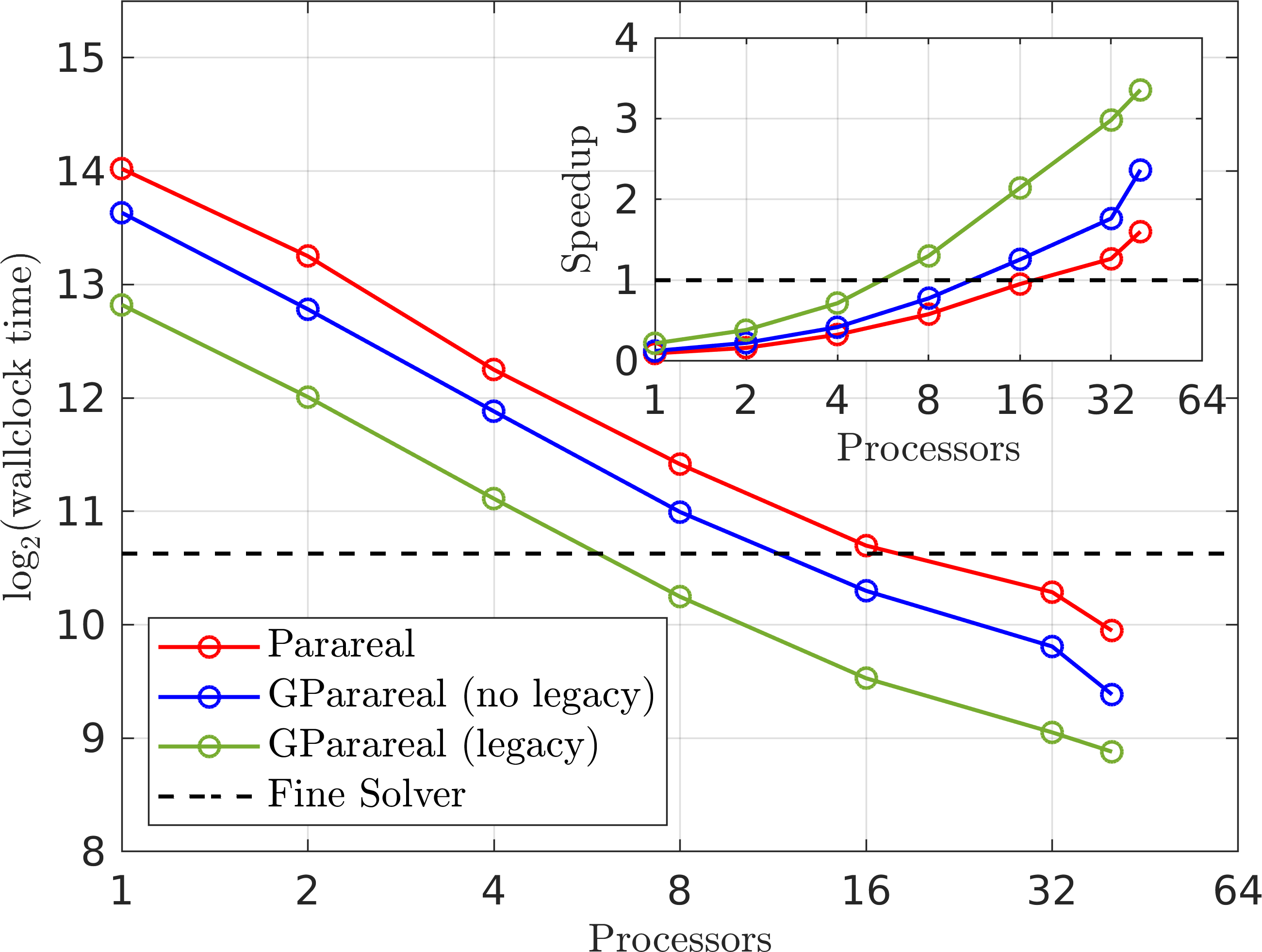}
        \caption{~}
    \end{subfigure}
    \caption{Numerical results obtained solving the R\"ossler system \eqref{eq:rossler} over $t \in [0,340]$.
    (a) Solutions from the fine solver, GParareal (with legacy data), and parareal plotted in the $(u_1,u_2)$-plane---both GParareal and parareal plotted only at times $\bm{t}$ for clarity.
    (b) The corresponding absolute errors between solutions from GParareal and parareal vs.\ the fine solution.
    (c) Maximum absolute errors \eqref{eq:stop} of each algorithm at successive iterations $k$ until tolerance $\varepsilon = 10^{-6}$ is met.
    (d) Median wallclock times (taken over 5 runs) of each simulation against the number of processors (up to 40). The inset plot displays the corresponding parallel speedup.}
    \label{fig:rossler}
\end{figure}
Next we solve the R{\"o}ssler system, 
\begin{equation} \label{eq:rossler}
    \frac{\rd u_1}{\rd t} = -u_2 - u_3, \quad \frac{\rd u_2}{\rd t} = u_1 + \hat{a} u_2, \quad \frac{\rd u_3}{\rd t} = \hat{b} + u_3(u_1 - \hat{c}),
\end{equation}
with parameters $(\hat{a},\hat{b},\hat{c}) = (0.2,0.2,5.7)$ that cause the system to exhibit chaotic behaviour \citep{rossler1976}.
We wish to integrate \eqref{eq:rossler} over $t \in [0,340]$ with initial values $\bm{u}_0 = (0,-6.78,0.02)^\intercal$ and solvers $\G=\text{RK1}$ and $\F=\text{RK4}$.
The interval is divided into $J=40$ time slices, $N_{\mathcal{G}}= 9 \times 10^{4}$ coarse steps, and $N_{\mathcal{F}}= 4.5 \times 10^{8}$ fine steps. The tolerance is set to $\varepsilon=10^{-6}$.

In this experiment, rather than obtaining legacy data by solving \eqref{eq:rossler} using alternative initial values (as we did in \cref{subsec:FHN}), we instead generate the data by integrating over a shorter time interval.
This is particularly useful if we are unsure how long to integrate our system for, i.e.\ to reach some long-time equilibrium state or reveal certain dynamics of the system, as is the case in many real-world dynamical systems. 
For example, many dynamical systems that feature random noise may exhibit metastability, in which trajectories spend (a long) time in certain states (regions of phase space) before transitioning to another state \citep{legoll2021a,grafke2017}.  
Such rare metastability may not be revealed/observed until the system has been evolved over a sufficiently large time interval.
We propose integrating over a `short' time interval, assessing the relevant characteristics of the solution obtained, and then integrating over a longer time interval (using the legacy data) if required.
Note that to do this, all parameters in both simulations must remain the same, with the exception of the time step widths---to ensure the legacy data is usable in the GP emulator in the longer simulation. 
Suppose we solve \eqref{eq:rossler} over $t \in [0,170]$, then we need to reduce $J$, $N_{\mathcal{G}}$, and $N_{\mathcal{F}}$ by a factor of two, i.e.\ use $J^{(2)}=J/2$, $N_{\mathcal{G}}^{(2)}=N_{\mathcal{G}}/2$, and $N_{\mathcal{F}}^{(2)}=N_{\mathcal{F}}/2$ in the shorter simulation. 

The legacy simulation, integrating over $[0,170]$, takes nine iterations to converge using GParareal (ten for parareal), giving us approximately $kJ^{(2)} = 9 \times 20 = 180$ legacy evaluations of $\F-\G$ (results not shown).
Integrating \eqref{eq:rossler} over the full interval $[0,340]$,  GParareal converges in four iterations sooner with the legacy data than without---refer to \cref{fig:rossler}(c).
In \cref{fig:rossler}(d) we can see that using the legacy data achieves a higher numerical speedup ($3.4\times$) compared to parareal ($1.6\times$).
In \cref{fig:rossler}(a) we see the trajectories from each simulation converging toward the R{\"o}ssler attractor and \cref{fig:rossler}(b) illustrates GParareal retaining a similar numerical accuracy to parareal with and without the legacy data.
Note the steadily increasing errors for both algorithms is due to the chaotic nature of the R{\"o}ssler system.

\subsection{Nonautonomous system} \label{subsec:nonauto}
\begin{figure}[b!]
    \centering
    \begin{subfigure}{0.49\linewidth}
        \includegraphics[width=\textwidth]{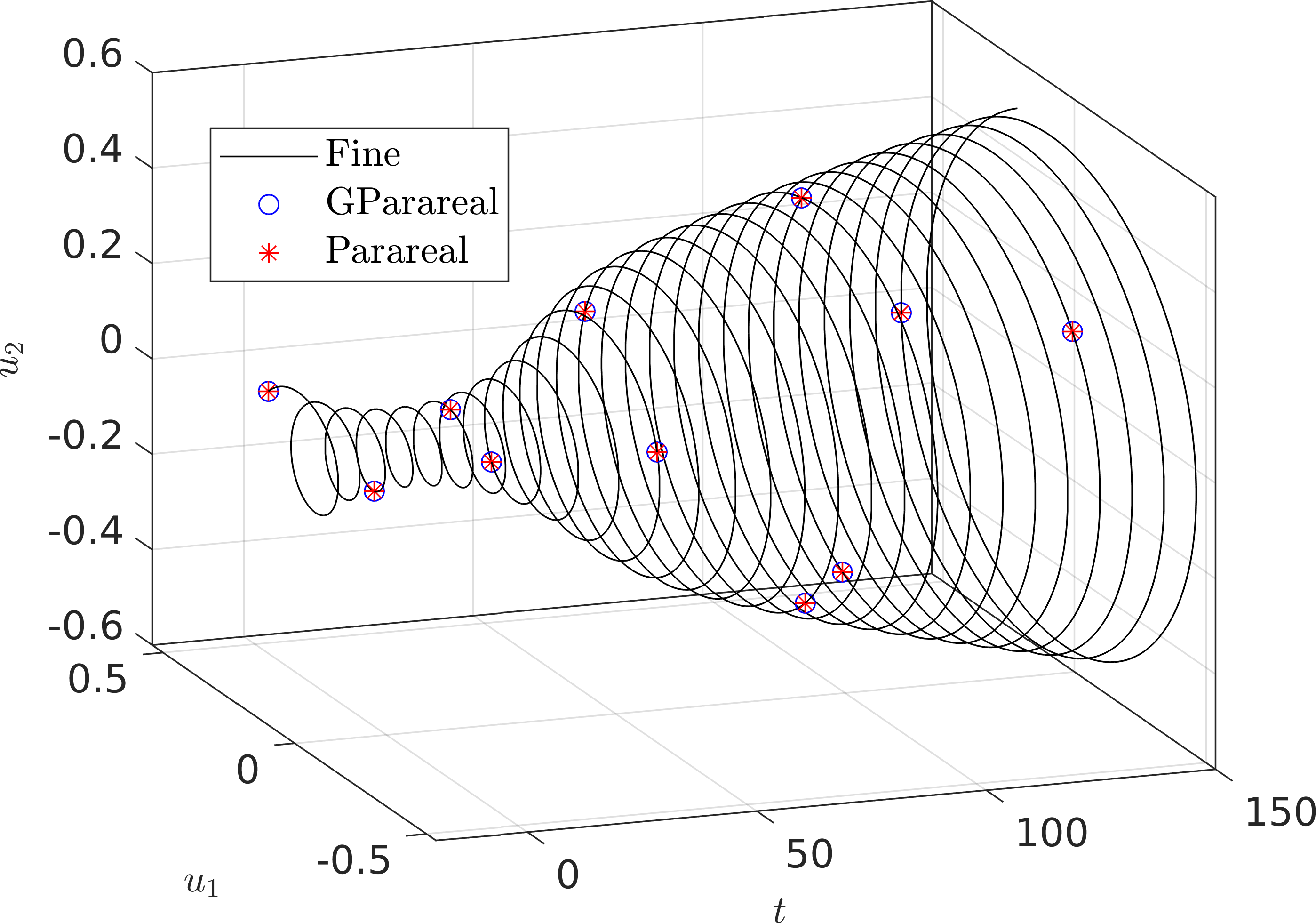}
        \caption{~}
    \end{subfigure}
    \begin{subfigure}{0.49\linewidth}
        \includegraphics[width=\textwidth]{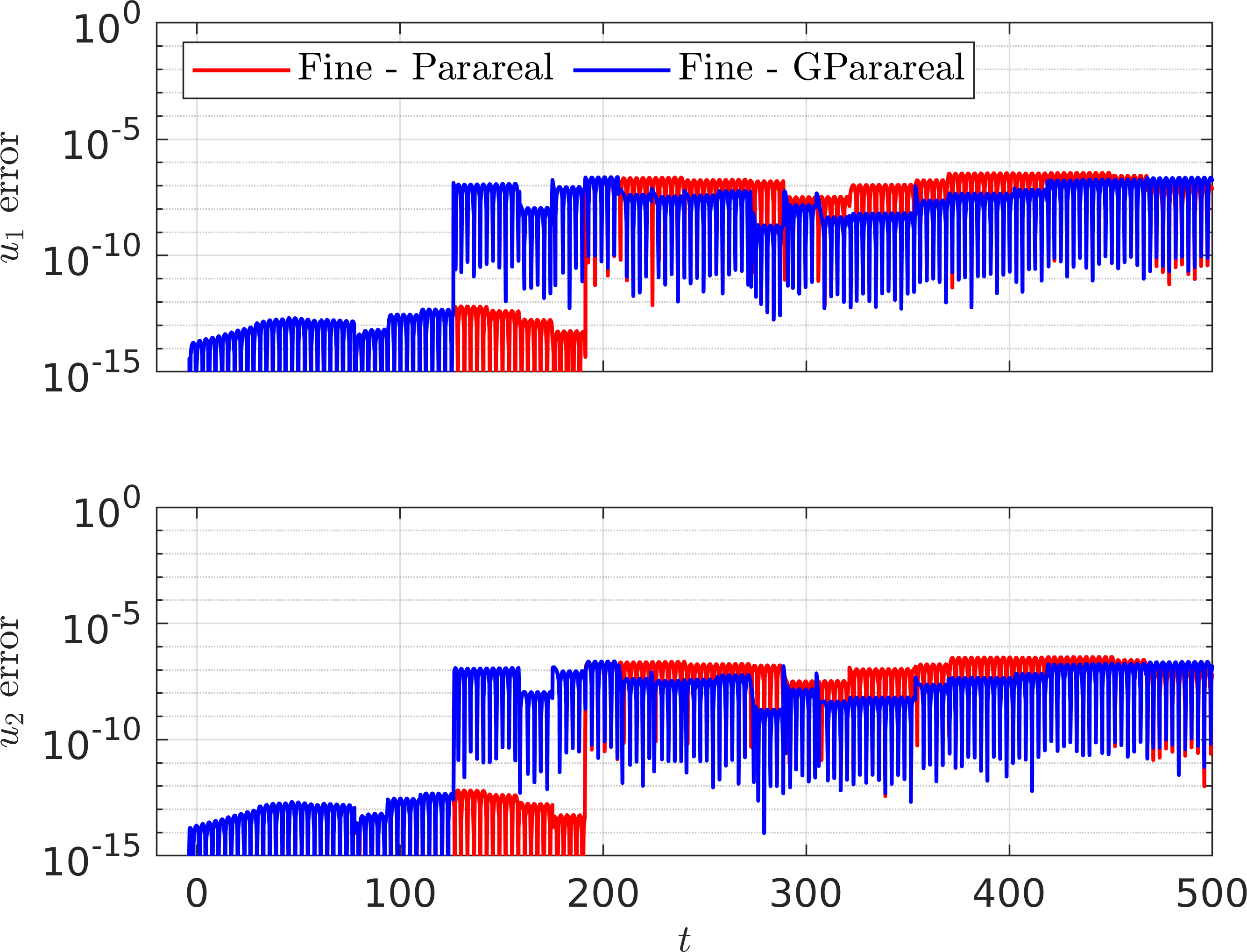}
        \caption{~}
    \end{subfigure}
    \caption{Numerical results obtained solving the nonautonomous system \eqref{eq:nonauto}.
    (a) Time-dependent solutions using the fine solver, GParareal, and parareal---both GParareal and parareal plotted only at times $\bm{t}$ on $[-20,150]$ for clarity. 
    (b) The corresponding absolute errors between solutions from GParareal and parareal vs.\ the fine solution, having converged after $10$ and $20$ iterations, respectively.}
    \label{fig:nonauto}
\end{figure}
Next, we consider a nonautonomous system of ODEs to demonstrate how GParareal handles explicit time dependence. 
We solve
\begin{equation} \label{eq:nonauto}
    \frac{\rd u_1}{\rd t} = -u_2 + u_1 \bigl(\frac{t}{500} - u_1^2 - u_2^2 \bigr), \quad \frac{\rd u_2}{\rd t} = u_1 + u_2 \bigl(\frac{t}{500} - u_1^2 - u_2^2 \bigr),
\end{equation}
over $t \in [-20,500]$---adapted from \cite{trefethen2017}.
As described in \cref{subsec:vectorised}, we transform this two-dimensional nonautonomous system into a three-dimensional autonomous system by introducing an additional variable $u_3(t) = t$, where $\rd u_3 / \rd t = 1$.
Given that we know $u_3(t)$ explicitly, the third dimension of $\F-\G$ need not be modelled with a GP. However, given the GPs are run in parallel anyway, this does not add to the cost of running GParareal.
\begin{table}[b!]
\caption{Wallclock time and speedup results obtained solving the nonautonomous system \eqref{eq:nonauto} for $J \in \{32,64,128,256,512\}$, where $k_{\text{para}}$ and $k_{\text{GPara}}$ are the number of iterations taken for parareal and GParareal to converge, respectively. (Top) Numerical results obtained upon simulation. (Bottom) Theoretical results calculated using \eqref{eq:Tpara}, \eqref{eq:TGPara}, \eqref{eq:Spara}, and \eqref{eq:SGPara}. All timings are measured in seconds.}
\label{table:1}
\begin{center}

\begin{tabular}{|c|c|c|c|c|c|c|c|c|c|c|} \hline \hline
$J$ & $k_{\text{para}}$ & $k_{\text{GPara}}$ & $T_{\mathcal{G}}$ & $T_{\mathcal{F}}$ & $T_{\text{GP}}$ & $T_{\text{serial}}$ & $T_{\text{para}}$ & $T_{\text{GPara}}$ & $S_{\text{para}}$ & $S_{\text{GPara}}$ \\ \hline
$32$ & $20$ & $10$ & $1.60\negE4$ & $4.23\posE3$ & $5.37$ & $1.35\posE5$ & $8.92\posE4$ & $4.33\posE4$ & $1.52$ & $3.13$ \\ \hline
$64$ & $31$ & $14$ & $9.80\negE5$ & $2.10\posE3$ & $18.36$ & $1.35\posE5$ & $6.75\posE4$  & $3.20\posE4$ & $2.00$ & $4.21$ \\ \hline
$128$ & $55$ & $16$ & $9.10\negE5$ & $1.06\posE3$ & $2.26\posE2$ & $1.35\posE5$ & $6.47\posE4$ & $1.90\posE4$ & $2.09$ & $7.13$ \\ \hline
$256$ & $99$ & $18$ & $6.90\negE5$ & $5.23\posE2$ & $1.02\posE2$ & $1.34\posE5$ & $5.64\posE4$  & $1.17\posE4$ & $2.37$ & $11.42$ \\ \hline
$512$ & $151$ & $15$ & $6.30\negE5$ & $2.62\posE2$ & $1.09\posE4$ & $1.34\posE5$ & $4.42\posE4$  & $2.10\posE4$ & $3.03$ & $6.39$ \\
\hline \hline
\end{tabular}
\bigskip

\begin{tabular}{|c|c|c|c|c|c|c|} \hline \hline
$J$ & $T_{\text{para}}$ & $T_{\text{GPara}}$ & $S_{\text{para}}$ & $S_{\text{GPara}}$ \\ \hline
$32$ & $8.47\posE4$ & $4.23\posE4$ & $1.60$ & $3.20$ \\ \hline
$64$ & $6.52\posE4$ & $2.95\posE4$ & $2.06$ & $4.57$ \\ \hline
$128$ & $5.81\posE4$ & $1.71\posE4$ & $2.33$ & $7.89$ \\ \hline
$256$ & $5.17\posE4$ & $9.51\posE3$ & $2.59$ & $14.07$ \\ \hline
$512$ & $3.95\posE4$ & $1.49\posE4$ & $3.39$ & $9.02$ \\
\hline \hline
\end{tabular}
\end{center}
\end{table}

We select solvers $\G=\text{RK1}$ and $\F=\text{RK8}$ with $N_{\mathcal{G}} = 2048$ and $N_{\mathcal{F}} = 5.12 \times 10^{5}$ steps, respectively.
We use $J=32$ time slices, initial condition $\bm{u}_0 = (0.1,0.1,-20)^{\intercal}$, and a stopping tolerance of $\varepsilon = 10^{-6}$.
In \cref{fig:nonauto}, we plot the solutions and corresponding errors generated by each of the solvers over time.
Again, the results illustrate good convergence to the fine solver solution, with GParareal taking $10$ iterations to locate the solution and parareal taking $20$.
We suspect that the superior performance of GParareal is partially due to the almost periodic nature of the solutions in \cref{fig:nonauto}(a), enabling the emulator to reproduce the dynamics of $\F-\G$ quite well.

\begin{figure}[t!]
    \centering
    \begin{subfigure}{0.49\linewidth}
        \includegraphics[width=\textwidth]{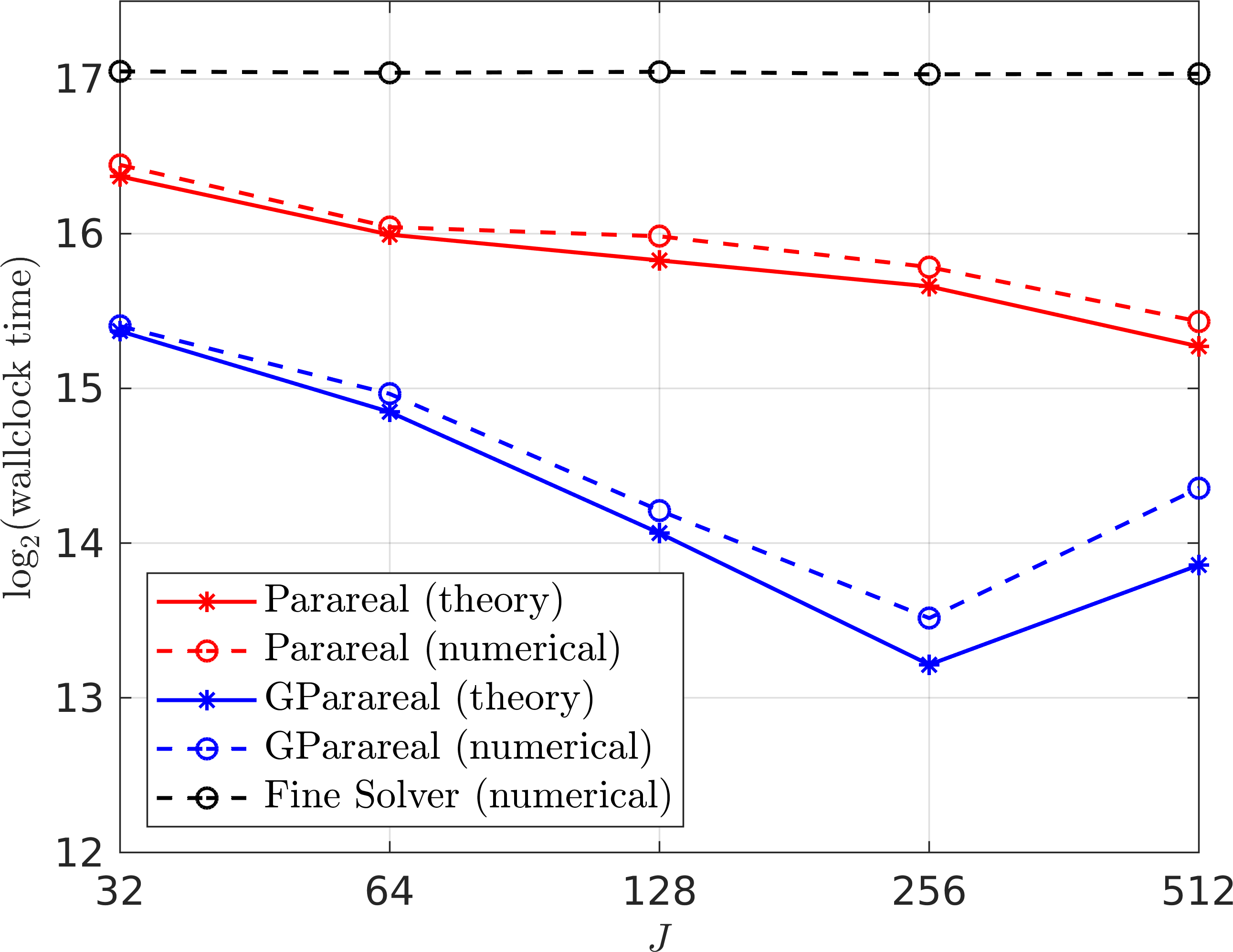}
        \caption{~}
    \end{subfigure}
    \begin{subfigure}{0.49\linewidth}
        \includegraphics[width=\textwidth]{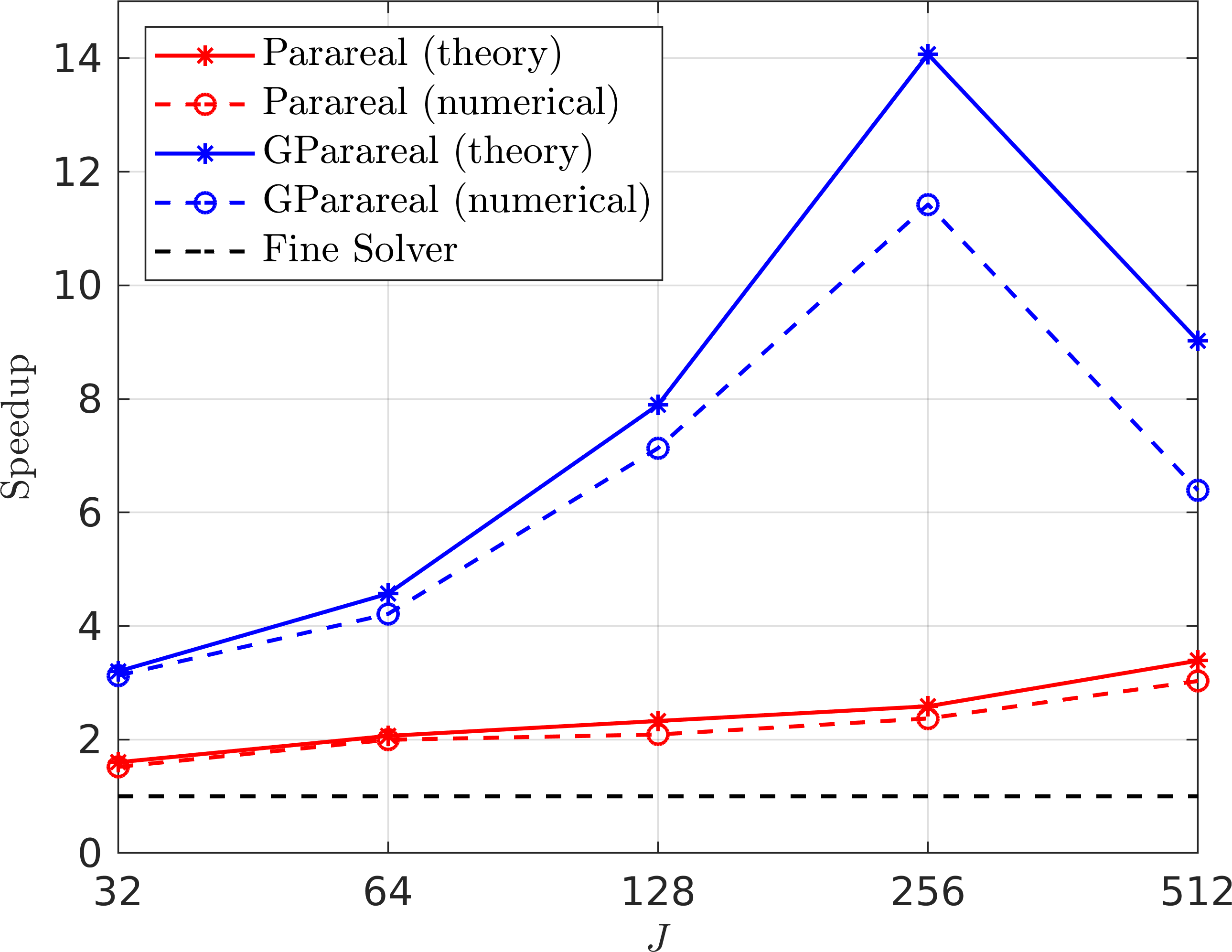}
        \caption{~}
    \end{subfigure}
    \caption{Numerical results obtained solving the nonautonomous system \eqref{eq:nonauto} for $J \in \{32,64,128,256,512\}$.
    (a) Wallclock times using the fine solver (dashed black), GParareal (dashed blue), and parareal (dashed red). Corresponding theoretical results are given by the solid lines, calculated using \eqref{eq:Tpara} and \eqref{eq:TGPara} for parareal and GParareal, respectively.  
    (b) The corresponding speedup results using the same lines and colours. The theoretical results were calculated using \eqref{eq:Spara} and \eqref{eq:SGPara} for parareal and GParareal, respectively.}
    \label{fig:nonauto2}
\end{figure}

Next, we run a series of simulations to measure the effect of increasing the number of time slices $J$ (and hence processors) on convergence, wallclock times, and speedup---see \cref{table:1}.
To do this, we increase the number of fine time steps to $N_{\mathcal{F}} = 5.12 \times 10^{10}$, so that $\F$ is sufficiently expensive to observe speedup.
We observe a good match between the numerical and theoretical results, presented in the top and bottom tables of \cref{table:1}, respectively, and visualised graphically in \cref{fig:nonauto2}.
Firstly, notice that $k_{\text{para}}$ increases with $J$ whilst $k_{\text{GPara}}$ remains largely unaffected, leading to speedups for GParareal being roughly $2\times$ to $4 \times$ that of parareal, up to $J = 256$. 
For both algorithms, the cost of $T_{\mathcal{G}}$ and $T_{\mathcal{F}}$ decreases as $J$ increases (due to fewer time steps per time slice), whilst $T_{\text{GP}}$ increases (due to increasing numbers of data points each simulation).
Note the exception of $T_{\text{GP}}=1.02\posE2$ when $J=256$ because hyperparameter optimisation converged within a few iterations and was therefore not carried out after this. 
Up to $J=256$, $T_{\text{GP}} < T_{\mathcal{F}}$ and so we observe increasing parallel speedup for GParareal.
When $J=512$, the cost of the GP overtakes that of $\F$ and so parallel speedup decreases, albeit still being double that of parareal.
Recall that if $T_{\text{GP}} > T_{\mathcal{F}}$, we may not opt to use GParareal in the first place, for the reasons outlined in \cref{subsec:complexity}.

\subsection{Double pendulum system} \label{subsec:double_pen}
\begin{figure}[t!]
    \centering
    \includegraphics[width=0.3\textwidth]{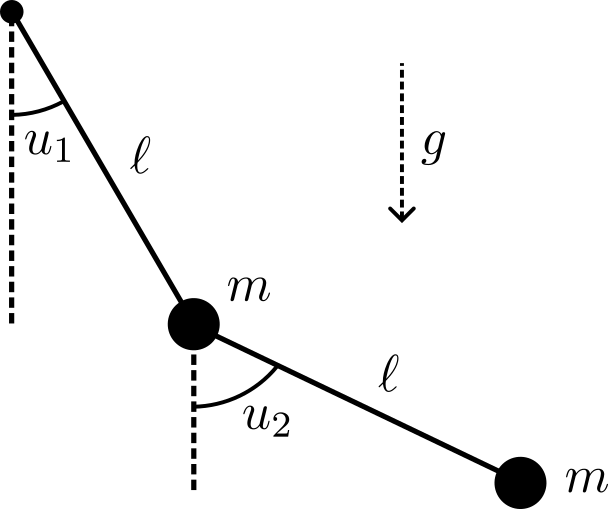}
    \caption{A schematic of the double pendulum system.}
    \label{fig:double_pen}
\end{figure}

Consider now the double pendulum system: a simple pendulum of mass $m$, rod length $\ell$, connected to another simple pendulum of equal mass $m$, rod length $\ell$, acting under gravity $g$ (see \cref{fig:double_pen}).
Four ODEs govern the dynamics of this system:
\begin{equation}
\begin{aligned} \label{eq:doub_pen}
    \frac{\rd u_1}{\rd t} &= u_3, \\
    \frac{\rd u_2}{\rd t} &= u_4, \\
    \frac{\rd u_3}{\rd t} &= \frac{-u_3^2 f_1(u_1,u_2) - u_4^2 \sin(u_1 - u_2) - 2 \sin(u_1) + \cos(u_1 - u_2) \sin(u_2)}{f_2(u_1,u_2)}, \\
    \frac{\rd u_4}{\rd t} &= \frac{2 u_3^2 \sin(u_1 - u_2) + u_4^2 f_1(u_1,u_2) + 2 \cos(u_1 - u_2) \sin(u_1) - 2 \sin(u_2)}{f_2(u_1,u_2)},
\end{aligned}
\end{equation}
where $f_1(u_1,u_2) = \sin(u_1 - u_2) \cos(u_1 - u_2)$ and $f_2(u_1,u_2) = 2 - \cos^2(u_1 - u_2)$ \citep{danby1997}.
Note that $m$, $\ell$, and $g$ have been scaled out of \eqref{eq:doub_pen} by letting $\ell = g$.
The variables $u_1$ and $u_2$ measure the angles between each pendulum and the vertical axis, while $u_3$ and $u_4$ measure the corresponding angular velocities.
\begin{figure}[b!]
    \centering
    \begin{subfigure}{0.48\linewidth}
        \includegraphics[width=\textwidth]{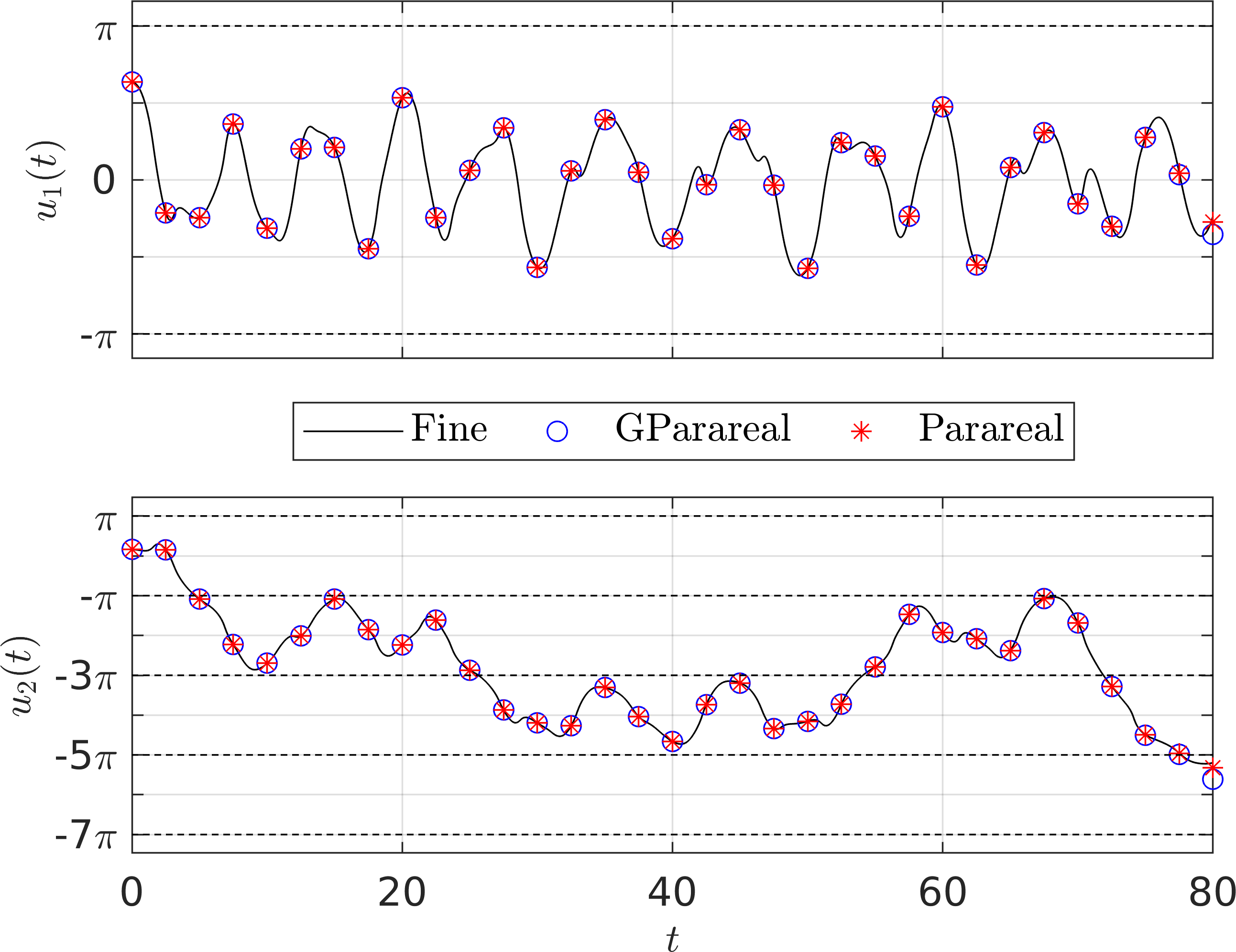}
        \caption{~}
    \end{subfigure}
    \begin{subfigure}{0.49\linewidth}
        \includegraphics[width=\textwidth]{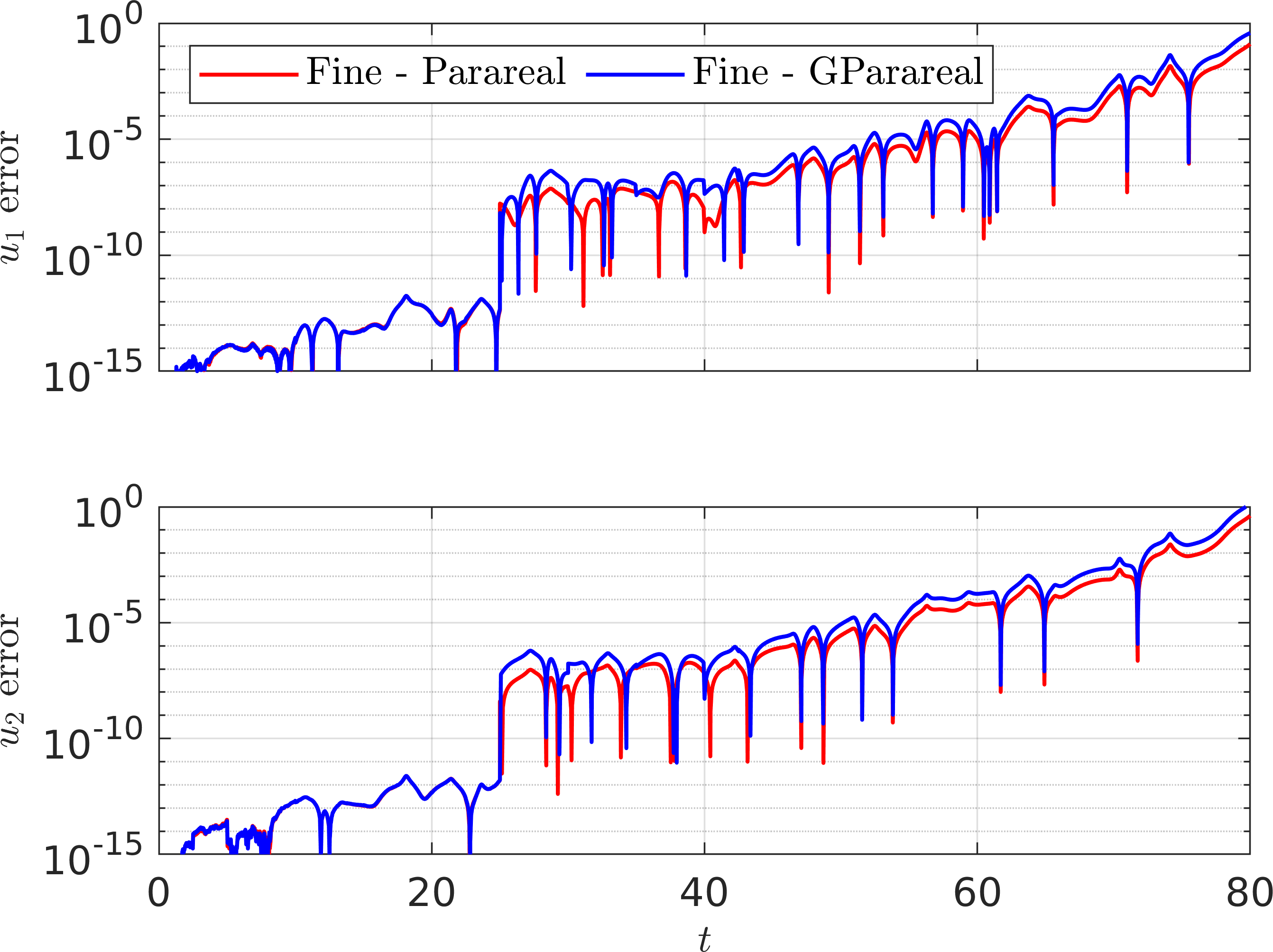}
        \caption{~}
    \end{subfigure}
    \caption{Numerical results obtained solving the double pendulum system \eqref{eq:doub_pen}.
    (a) Time-dependent solutions for $u_1$ and $u_2$ using the fine solver, GParareal, and parareal---both GParareal and parareal plotted only at times $\bm{t}$ for clarity. Dashed lines indicate ``turning over'' angles, at which either pendulum passes through an odd multiple of $\pi$. 
    (b) The corresponding absolute errors between solutions from GParareal and parareal vs.\ the fine solution, having converged after $23$ and $22$ iterations, respectively.}
    \label{fig:doubpen2}
\end{figure}

\begin{figure}[t!]
    \centering
    \begin{subfigure}{0.49\linewidth}
        \includegraphics[width=\textwidth]{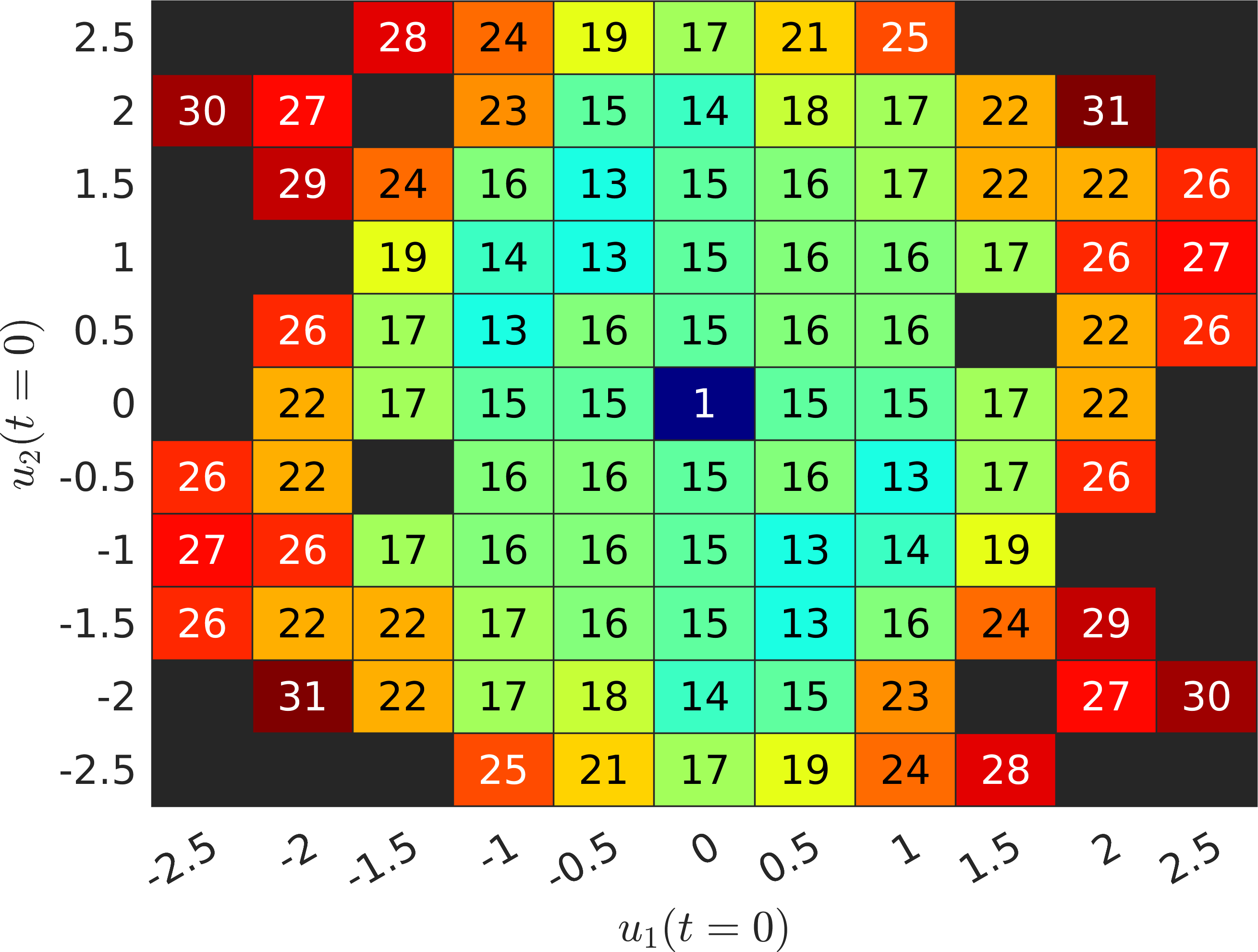}
        \caption{Parareal}
    \end{subfigure}
    \begin{subfigure}{0.49\linewidth}
        \includegraphics[width=\textwidth]{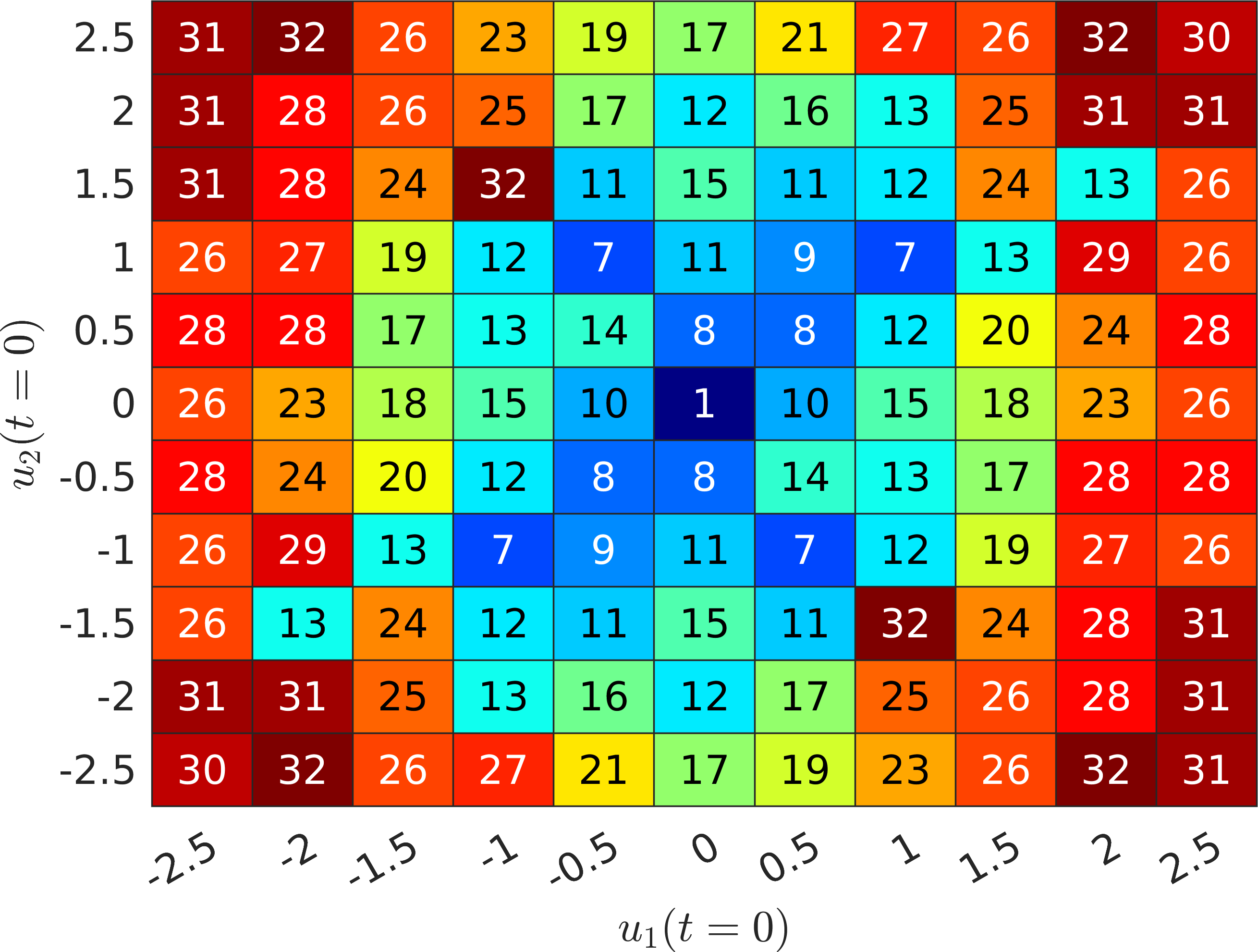}
        \caption{GParareal}
    \end{subfigure}
    \caption{Heat maps displaying the number of iterations taken until convergence $k$ of (a) parareal and (b) GParareal when solving the double pendulum system \eqref{eq:doub_pen} for different initial angles $(u_1(0),u_2(0)) \in [-2.5,2.5]^2$ and initial angular velocities $(u_3(0),u_4(0)) = (0,0)$, i.e.\ the pendulums are released from rest. 
    Black boxes indicate where parareal returned a \texttt{NaN} value during simulation.}
    \label{fig:doub_pen_heat}
\end{figure}

For this experiment, we select solvers $\G=\text{RK1}$ and $\F=\text{RK8}$ with $N_{\mathcal{G}} = 3072$ and $N_{\mathcal{F}} = 2.1504 \times 10^5$ steps, respectively.
We integrate over $t \in [0,80]$, using $J=32$ time slices with a stopping tolerance $\varepsilon = 10^{-6}$.
In \cref{fig:doubpen2}, we plot solutions for $u_1$ and $u_2$ over time using initial conditions $\bm{u}_0 = (2,0.5,0,0)^{\intercal}$, i.e. the pendulums are positioned at some (positive) initial angles and released from rest.
Observe how both pendulums move chaotically, with the inner pendulum oscillating within $[-\pi,\pi]$ and the outer pendulum oscillating between odd multiples of $\pi$, ``turning over'' a number of times. 
We attain good solution accuracy from GParareal with respect to the fine solution with errors slowly increasing over time due to the chaotic nature of the system, much like what was seen in the R{\"o}ssler experiments in \cref{subsec:rossler}.
We plot $k$ for various initial angles in \cref{fig:doub_pen_heat} to highlight the system's sensitivity to initial conditions.
For small initial angles, GParareal converges sooner than parareal, but for much larger angles both algorithms use almost all of the $32$ iterations to locate a solution (and in some cases, parareal does not return a solution).

\begin{figure}[t!]
    \centering
    \begin{subfigure}{0.49\linewidth}
        \includegraphics[width=\textwidth]{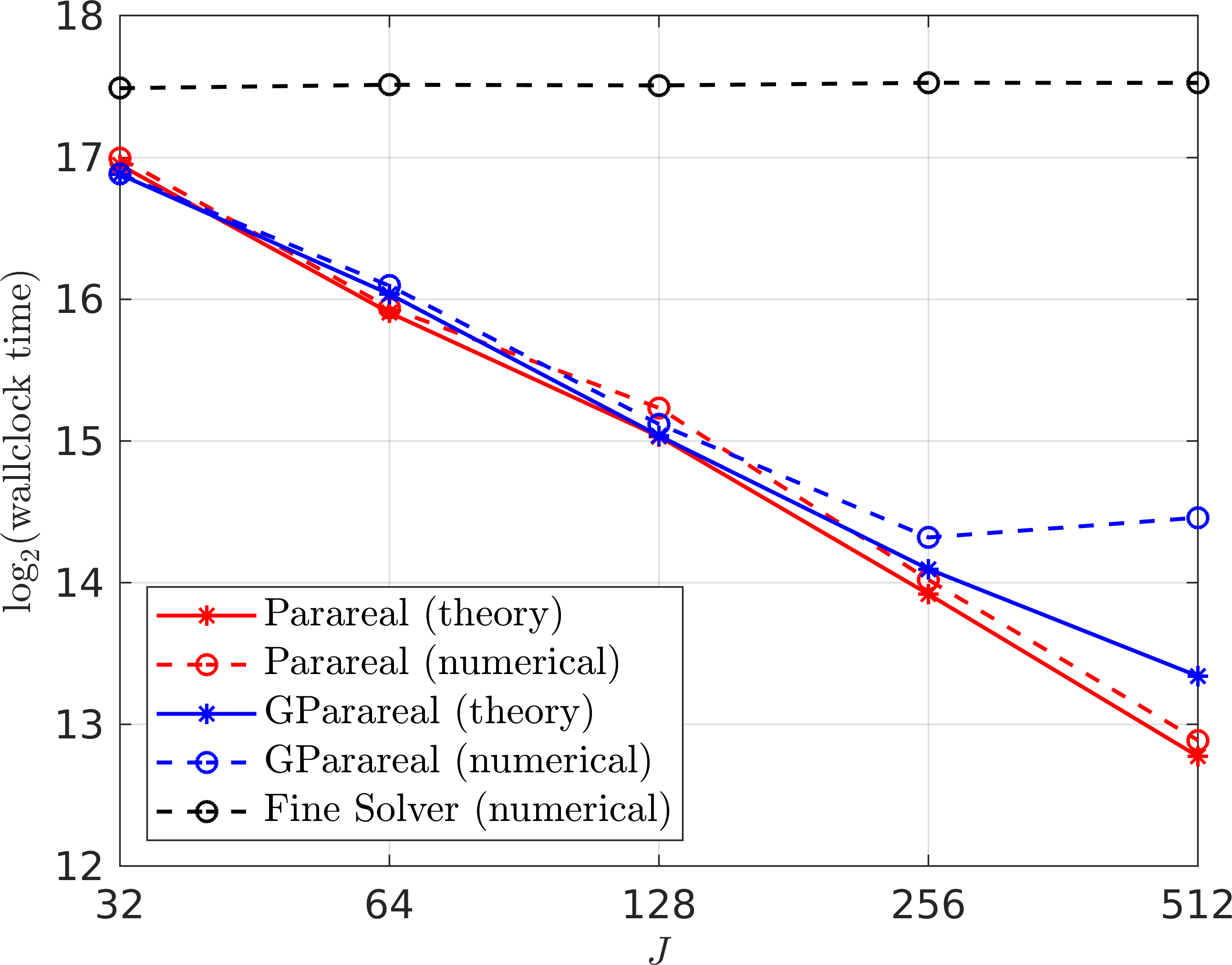}
        \caption{~}
    \end{subfigure}
    \begin{subfigure}{0.49\linewidth}
        \includegraphics[width=\textwidth]{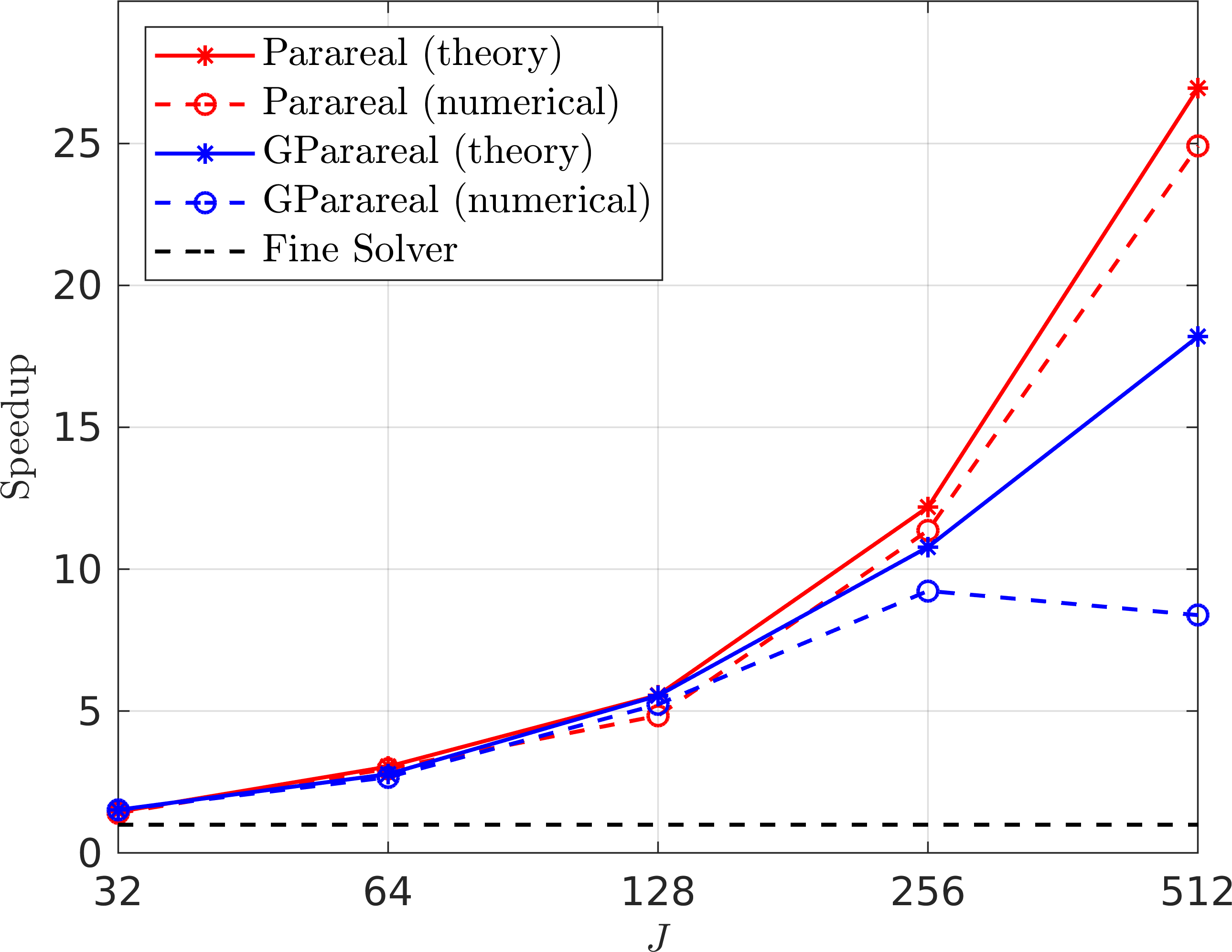}
        \caption{~}
    \end{subfigure}
    \caption{Numerical results obtained solving the double pendulum system \eqref{eq:doub_pen} for $J \in \{32,64,128,256,512\}$.
    (a) Wallclock times using the fine solver (dashed black), GParareal (dashed blue), and parareal (dashed red). Corresponding theoretical results are given by the solid lines, calculated using \eqref{eq:Tpara} and \eqref{eq:TGPara} for parareal and GParareal, respectively.  
    (b) The corresponding speedup results using the same lines and colours. The theoretical results were calculated using \eqref{eq:Spara} and \eqref{eq:SGPara} for parareal and GParareal, respectively.}
    \label{fig:doubpen3}
\end{figure}

In \cref{table:2} and \cref{fig:doubpen3}, we again test the effect of increasing $J$ on wallclock times, speedup, and convergence. To do this, we increase the number of fine time steps to $N_{\mathcal{F}} = 2.1504 \times 10^{10}$.
We purposefully choose an initial condition ($\bm{u}_0$ above) for which both algorithms converge in approximately the same number of iterations, so that we can directly observe how the increasing GP cost affects the performance of GParareal for large $J$.
Under these circumstances, we can think of the wallclock time for GParareal as (approximately) the wallclock time of parareal plus the wallclock time of the GP conditioning/optimisation. 
For $J \leq 128$, we observe that $T_{\text{GP}} < T_{\mathcal{F}}$ and so the speedup of GParareal and parareal are approximately the same.
In these cases, using GParareal is no more costly than using parareal, with the additional benefit of being able to re-use the acquisition data for a future simulation, if needed. 
For $J \geq 256$, we begin to observe $T_{\text{GP}} \approx T_{\mathcal{F}}$ (or larger), so the numerical speedup of GParareal begins to plateau.
We should reiterate, however, that using so many processors for such a small test problem is quite excessive.

\begin{table}[b!]
\caption{Wallclock time and speedup results obtained solving the double pendulum system \eqref{eq:doub_pen} for $J \in \{32,64,128,256,512\}$, where $k_{\text{para}}$ and $k_{\text{GPara}}$ are the number of iterations taken for parareal and GParareal to converge, respectively. (Top) Numerical results obtained upon simulation. (Bottom) Theoretical results calculated using \eqref{eq:Tpara}, \eqref{eq:TGPara}, \eqref{eq:Spara}, and \eqref{eq:SGPara}. All timings are measured in seconds.}
\label{table:2}
\begin{center}

\begin{tabular}{|c|c|c|c|c|c|c|c|c|c|c|} \hline \hline
$J$ & $k_{\text{para}}$ & $k_{\text{GPara}}$ & $T_{\mathcal{G}}$ & $T_{\mathcal{F}}$ & $T_{\text{GP}}$ & $T_{\text{serial}}$ & $T_{\text{para}}$ & $T_{\text{GPara}}$ & $S_{\text{para}}$ & $S_{\text{GPara}}$ \\ \hline
$32$ & $22$ & $21$ & $2.53\negE4$ & $5.75\posE3$ & $19.75$ & $1.84\posE5$ & $1.31\posE5$ & $1.21\posE5$ & $1.41$ & $1.52$ \\ \hline
$64$ & $21$ & $23$ & $1.46\negE4$ & $2.93\posE3$ & $17.66$ & $1.87\posE5$ & $6.29\posE4$ & $7.00\posE4$ & $2.97$ & $2.67$ \\ \hline
$128$ & $23$ & $23$ & $1.27\negE4$ & $1.46\posE3$ & $1.20\posE2$ & $1.86\posE5$ & $3.85\posE4$ & $3.56\posE4$ & $4.84$ & $5.24$ \\ \hline
$256$ & $21$ & $23$ & $9.10\negE5$ & $7.35\posE2$ & $5.45\posE2$ & $1.89\posE5$ & $1.66\posE4$ & $2.04\posE4$ & $11.36$ & $9.24$ \\ \hline
$512$ & $19$ & $22$ & $7.00\negE5$ & $3.69\posE2$ & $2.26\posE3$ & $1.89\posE5$ & $7.58\posE3$ & $2.25\posE4$ & $24.90$ & $8.38$ \\
\hline \hline
\end{tabular}
\bigskip

\begin{tabular}{|c|c|c|c|c|c|c|} \hline \hline
$J$ & $T_{\text{para}}$ & $T_{\text{GPara}}$ & $S_{\text{para}}$ & $S_{\text{GPara}}$ \\ \hline
$32$ & $1.26\posE5$ & $1.21\posE5$ & $1.45$ & $1.52$ \\ \hline
$64$ & $6.14\posE4$ & $6.72\posE4$ & $3.05$ & $2.78$ \\ \hline
$128$ & $3.35\posE4$ & $3.36\posE4$ & $5.57$ & $5.54$ \\ \hline
$256$ & $1.55\posE4$ & $1.75\posE4$ & $12.19$ & $10.78$ \\ \hline
$512$ & $7.01\posE3$ & $1.04\posE4$ & $26.94$ & $18.20$ \\
\hline \hline
\end{tabular}
\end{center}

\end{table}

\newpage
\section{Discussion} \label{sec:discussion}
In this paper, we present a time-parallel algorithm (GParareal) that iteratively locates a numerical solution to a system of ODEs. 
It does so using a predictor-corrector, comprised of numerical solutions from coarse ($\G$) and fine ($\F$) integrators. 
However, unlike the classical parareal algorithm, it uses a Gaussian process (GP) emulator to infer the correction term $\F-\G$.
The numerical experiments reported in \cref{sec:numerics} demonstrate that GParareal performs favourably compared to parareal, converging in fewer iterations and achieving increased parallel speedup for a number of low-dimensional ODE systems.
We also demonstrate how GParareal can make use of legacy data, i.e.\ prior $\F$ and $\G$ data obtained during a previous simulation of the same system (using different ICs or a shorter time interval), to pre-train the emulator and converge even faster---something that existing time-parallel methods cannot do. 

In \cref{subsec:FHN}, using just the data obtained during simulation (acquisition data), GParareal achieves an almost two-fold increase in speedup over parareal when solving the FitzHugh-Nagumo model.
Simulating over a range of initial values, GParareal converged in fewer than half the iterations taken by parareal and, in some cases, managed to converge when the coarse solver was too poor for parareal.
When using legacy data, GParareal could converge in even fewer iterations.
Similar results were illustrated for the R{\"o}ssler system in \cref{subsec:rossler} but with legacy data obtained from a prior simulation over a shorter time interval---beneficial when one does not know how long to integrate a system for.
In \cref{subsec:nonauto,subsec:double_pen}, GParareal was tested on a larger number of processors (up to $512$), verifying the theoretical computational complexity results given in \cref{subsec:complexity} and that the cost of the GP needs to be much smaller than the cost of the fine solver in order for speedup to be maximised.
In all cases, the solutions generated by GParareal were of a numerical accuracy comparable to those found using parareal. 

In its current implementation, GParareal may, however, suffer from the curse of dimensionality in two ways.
First, an increasing number of data points, $\mathcal{O}(kJ)$, is problematic for the standard cubic complexity GP implemented here. 
In this case, a more sophisticated (non-cubic complexity) emulator or perhaps using neural networks could be beneficial.  
Second, trying to emulate a $d$-dimensional function $\F-\G$ is difficult if the number of evaluation points is not sufficient.
One option to tackle this may be to obtain more acquisition data by launching more $\F$ and $\G$ runs using the idle processors to further train the emulator at little additional computational cost.
However, as shown in \cref{subsec:convergence}, the accuracy of the GP emulator is strongly controlled by the fill distance of the set of evaluation points, which is generally difficult to restrict when $d$ is large.
One could think about using legacy data generated by evaluating $\F-\G$ at specific input locations (for example, a uniform grid) that satisfy certain fill distance requirements in the state space. 

It should also be noted that GParareal may not always provide faster convergence using legacy data if such legacy evaluations of $\F-\G$ lay `far away', i.e.\ over one or two input length-scales away, from the initial values of interest in the current simulation.
In this case, GParareal would rely more heavily on its acquisition data.
There is no immediate remedy for such a situation, but using a fallback parareal correction, as suggested in the next paragraph, could be an option. 

In equation \eqref{eq:V2}, we approximate a Gaussian distribution by taking its expected value, ignoring uncertainty in the GP posterior for $\F-\G$.
In this setting, the GP emulator is used to interpolate the $\F-\G$ data, hence it is perfectly acceptable to swap it out for any other sufficiently accurate interpolation method, e.g.\ kernel ridge regression \citep{kanagawa2018}.
During early iterations of GParareal, when little acquisition data are available, the uncertainty in the GP posterior (i.e\ the variance) may be large at points of interest.
By retaining the GP posterior uncertainty, one could (ideally) propagate the full uncertainty using the coarse solver to the next time step and then continue. 
While this would produce a probabilistic version of GParareal, this would be a computationally expensive process that we wish to avoid at this stage. 
One alternative to approximating \eqref{eq:V2} by its expected value could be to draw a random sample instead.
A sampling-based solver such as this would return a stochastic solution to the ODE, much like the stochastic parareal algorithm presented in \cite{pentland2022}.
It is unclear how this algorithm would perform vs.\ parareal (or even stochastic parareal), however, it could still make use of legacy data following successive independent simulations. 
Another possible alternative to approximating \eqref{eq:V2} arises if the input initial value is at least one or two length-scale distances away from any other known input value in our acquisition dataset.
In this case, we then might expect the GP emulation of the mean in \eqref{eq:V2} to have high variance and so a fallback to the deterministic parareal correction for $\F-\G$ (see \eqref{eq:pred-corr}) could be built in as a next best correction to the coarse prediction in \eqref{eq:V2}.
Among others, these are two alternative formulations of GParareal that are worth investigating to account for the whole Gaussian distribution provided by the emulator and not just its mean value.

Follow-up work will focus on extending GParareal, using some of the methods suggested above, to solve higher-dimensional systems of ODEs in parallel (possibly PDEs).
In the longer term, we aim to develop a truly probabilistic time-parallel numerical method that can account for the inherent uncertainty in the GP emulator, returning a probability distribution rather than point estimates over the solution.


\section*{Acknowledgements}
KP is funded by the Engineering and Physical Sciences Research Council through the MathSys~II CDT (grant EP/S022244/1) as well as the Culham Centre for Fusion Energy.
TJS is partially supported by the Deutsche Forschungs\-gemeinschaft through project 415980428.
This work has partly been carried out within the framework of the EUROfusion Consortium and has received funding from the Euratom research and training programme 2014--2018 and 2019--2020 under grant agreement No.~633053.
The authors would also like to acknowledge the University of Warwick Scientific Computing Research Technology Platform for assistance in the research described in this paper, in particular Arkadiy Davydov.
The views and opinions expressed herein do not necessarily reflect those of any of the above-named institutions or funding agencies.
For the purpose of open access, the author has applied a CC BY public copyright licence to any
Author Accepted Manuscript version arising.


\bibliographystyle{abbrvnat}  
\bibliography{references}  


\newpage

\begin{appendices} \label{appendix}
\crefalias{section}{appendix}

\section{Psuedocode for GParareal}
\label{algorithm}

\SetKwComment{Comment}{\%}{}
\SetKwInput{KwInput}{Initialise}

\begin{algorithm}[H]
\caption{GParareal} 
\KwInput{Set counters $k = I = 0$ and define $V_j^k$, $\hat{V}_j^k$ and $\tilde{V}_j^k$ as the refined, coarse, and fine solutions at the $j$th mesh point and $k$th iteration respectively (note $V_0^k = \hat{V}_0^k = \tilde{V}_0^k = u^0 \ \forall k$). If known, initialise any legacy $\F-\G$ input data $\bm{x}$, output data $\bm{y}$, and hyperparameters $\bm{\theta}$.}

\Comment{Calculate approximate initial values at each $t_j$ by running $\G$ serially.}
\For{$j = 1$ \KwTo $J$}{
    $\hat{V}_j^0 = \G(\hat{V}_{j-1}^0)$\;
    $V_j^0 = \hat{V}_j^0$\;
}

\For{$k = 1$ \KwTo $J$}{
    \Comment{Propagate refined solutions (from iteration $k-1$) on unconverged sub-intervals by running $\F$ in parallel.}
    \For{$j = I+1$ \KwTo $J$}{
    $\tilde{V}_j^{k-1} = \F(V_{j-1}^{k-1})$\;
    }
    $I = I + 1$\;
    $V^k_I = \tilde{V}^{k-1}_I$ for all $k$ \Comment*[r]{copy converged solution at $t_I$ to future $k$.}
    $\bm{x} = \text{append}(\bm{x},(V^{k-1}_I,\dots,V^{k-1}_{J-1})^{\text{T}})$ \Comment*[r]{collect new input data.}
    $\bm{y} = \text{append}(\bm{y},(\tilde{V}^{k-1}_{I+1} - \hat{V}^{k-1}_{I+1},\dots,\tilde{V}^{k-1}_{J}-\hat{V}^{k-1}_{J})^{\text{T}})$ \Comment*[r]{collect new output data.}
    $\bm{\theta} = \text{GPoptimise}(\bm{x},\bm{y},\bm{\theta})$ \Comment*[r]{optimise hyperparameters.}
    \Comment{Propagate refined solution (at iteration $k$) with $\G$, then correct using the expected value of the GP prediction \eqref{eq:GP_post} (this step \emph{cannot} be carried out in parallel).}

\For{$j = I+1$ \KwTo $J$}{
    $x^{\star} = V_{j-1}^{k}$\;
    $\hat{V}_j^k = \G(x^{\star})$\; 
    $y^{\star} = \text{GPpredict}(\bm{x},\bm{y},\bm{\theta},x^{\star})$ \Comment*[r]{returns Gaussian random variable} 
    $V_{j}^{k} = \mathbb{E}[y^{\star}] + \hat{V}_{j}^{k}$\;
}
\Comment{Evaluate the stopping criterion, saving all solutions up to $t_I$. }
$I = \underset{n \in \{I,\dots,N\}}{\mathrm{max}}  \ | V^k_i - V^{k-1}_i | < \varepsilon \quad \forall i < n$\;
\If{$I = N$}{
    \textbf{return} $k$, $\bm{V}^k$, $\bm{x}$, $\bm{y}$, $\bm{\theta}$ \Comment*[r]{if tolerance met for all time steps, stop.}
}
}
\end{algorithm}


\end{appendices}


\end{document}